\documentclass[bpam,reqno]{ipart}

\RequirePackage{hyperref}

\startlocaldefs
\theoremstyle{plain}

\usepackage{graphicx}
\usepackage{subfigure}
\usepackage{amssymb}
\usepackage{amsmath}
\usepackage{amsfonts}
\usepackage{theorem}
\usepackage{mathrsfs}
\usepackage{mathtools}
\usepackage{bm}
\usepackage{color}
\usepackage{exscale}
\usepackage{relsize}
\usepackage{epstopdf}
\usepackage{float}
\usepackage{nicefrac}
\usepackage{booktabs,multirow} % for much better looking tables
\usepackage{array} % for better arrays (eg matrices) in maths
\usepackage{paralist} % very flexible & customisable lists (eg. enumerate/itemize, etc.)
\usepackage{verbatim} % adds environment for commenting out blocks of text & for better verbatim
\usepackage{subfigure} % make it possible to include more than one captioned figure/table in a single float

\usepackage{bm}

\usepackage{tikz}
\usetikzlibrary{positioning}
\usepackage{xcolor}

\allowdisplaybreaks[1]

\newcommand\eref[1]{(\ref{#1})}

\newcommand*\xbar[1]{%
  \hbox{%
    \vbox{%
      \hrule height 0.5pt % The actual bar
      \kern0.4ex%         % Distance between bar and symbol
      \hbox{%
        \kern-0.05em%      % Shortening on the left side
        \ensuremath{#1}%
        \kern-0.00em%      % Shortening on the right side
      }%
    }%
  }%
}

\newcommand{\mF}{\bm{F}}
\newcommand{\mK}{\bm{K}}

\newcommand{\mU}{\bm{U}}

\newcommand{\dt}{\Delta t}
\newcommand{\dx}{\Delta x}

\newcommand{\hf}{{\frac{1}{2}}}

\newcommand{\jph}{{j+\frac{1}{2}}}
\newcommand{\jmh}{{j-\frac{1}{2}}}

\def\softd{{\leavevmode\setbox1=\hbox{d}%
          \hbox to 1.05\wd1{d\kern-0.4ex{\char039}\hss}}}

\endlocaldefs

\pubyear{2024}
\volume{0}
\issue{0}
\firstpage{1}
\lastpage{28}
\begin{document}

\begin{frontmatter}
\title[WB A-WENO Schemes]{A Well-Balanced Fifth-Order A-WENO Scheme Based on Flux Globalization}

\begin{aug}
\author{\fnms{Shaoshuai} \snm{Chu}\ead[label=e1]{chuss2019@mail.sustech.edu.cn}},
\address{Department of Mathematics, RWTH Aachen University, Aachen, 52056\\
Germany\\
\printead{e1}}
\author{\fnms{Alexander} \snm{Kurganov}\ead[label=e2]{alexander@sustech.edu.cn}},
\address{Department of Mathematics, Shenzhen International Center for Mathematics, and Guangdong Provincial Key Laboratory of Computational
Science and Material Design, Southern University of Science and Technology, Shenzhen, 518055\\
China\\
\printead{e2}}
\and
\author{\fnms{Ruixiao} \snm{Xin}\ead[label=e3]{12331009@mail.sustech.edu.cn}}
\address{Department of Mathematics, Southern University of Science and Technology, Shenzhen, 518055\\
China\\
\printead{e3}}
\end{aug}
\received{\sday{1} \smonth{1} \syear{2024}}

\begin{abstract}
We construct a new fifth-order flux globalization based well-balanced (WB) alternative weighted essentially non-oscillatory (A-WENO) scheme for general nonconservative systems. The proposed scheme is a higher-order extension of the WB path-conservative central-upwind (PCCU) scheme recently proposed in [{\sc A. Kurganov, Y. Liu and R. Xin}, J. Comput. Phys., 474 (2023), Paper No. 111773]. We apply the new scheme to the nozzle flow system and the two-layer shallow water equations. We conduct a series of numerical experiments, which clearly demonstrate the advantages of using the fifth-order extension of the flux globalization based WB PCCU scheme.
\end{abstract}


%%  Upper case for every keyword

%\tableofcontents
\end{frontmatter}

\section{Introduction}
This paper is focused on the development of high-order well-balanced (WB) numerical methods for the nonconservative hyperbolic systems,
which in the one-dimensional (1-D) case read as
\begin{equation}
\bm U_t+\bm F(\bm U)_x=B(\bm U)\bm U_x.
\label{1.1}
\end{equation}
Here, $x$ is the spatial variable, $t$ is time, $\mU\in\mathbb R^d$ is a vector of unknowns, $\mF:\mathbb R^d\to\mathbb R^d$ is the
nonlinear flux function, and $B\in\mathbb R^{d\times d}$.

Developing a highly accurate numerical scheme for \eref{1.1} is a challenging task due to the presence of the nonconservative product terms
on the right-hand side (RHS) of the studied system, whose weak solutions can be understood as the Borel measures (see
\cite{DLM,LPG2002,LPG2004}) rather than in the sense of distributions. In addition, a good scheme for \eref{1.1} should be WB, that is, it
should be capable of exactly preserving some of the physically relevant steady-state solutions of (1.1). As it was shown in \cite{KLX_21},
for many nonconservative systems steady-state solutions satisfy
\begin{equation}
\bm F(\bm U)_x-B(\bm U)\bm U_x=M(\bm U)\bm E(\bm U)_x=\bm0,
\label{1.4}
\end{equation}
where $M\in\mathbb R^{d\times d}$ and $\bm E$ is the vector of equilibrium variables, which are constant at the steady states.

A popular approach for developing accurate and WB numerical schemes for \eref{1.1} is based on the concept of weak Borel measure solutions,
which leads to a class of path-conservative schemes; see, e.g., \cite{BDGI,CMP2017,Cha,CKM,PC2009,PGCCMP,SGBNP} and references therein. One
of these schemes, a second-order path-conservative central-upwind (PCCU) scheme \cite{CKM}, is particularly simple as it is based on
Riemann-problem-solver-free central-upwind (CU) schemes, which were originally developed in \cite{KLin,KNP,KTcl} as a ``black-box'' solver
for general multidimensional hyperbolic systems. One of the drawbacks of the PCCU schemes is that they can preserve simple steady states
only. For instance, in the two-layer shallow water equations, the PCCU schemes can only preserve the so-called ``lake-at-rest'' steady
states.

In \cite{KLX_21}, path-conservative techniques were incorporated into the flux globalization framework, which led to a flux globalization
based WB PCCU scheme, which is capable of preserving a much wider range of steady states while still accurately treating the nonconservative
products. In general, a flux globalization approach, which was introduced in \cite{CHO_18,CDH2009,DDMA11,GC2001,MGD11}, relies on the
following quasi-conservative form of \eref{1.1}:
\begin{equation}
\bm U_t+\bm K(\bm U)_x=\bm0,
\label{1.6}
\end{equation}
where $\bm K(\bm U)$ is a global flux
\begin{equation}
\bm K(\bm U)=\bm F(\bm U)-\bm R(\bm U),\quad\bm R(\bm U)=\int\limits_{\widehat x}^x\left[B(\bm U(\xi,t))\bm U_\xi(\xi,t)\right]\,{\rm d}\xi,
\label{1.7}
\end{equation}
and ${\widehat x}$ is an arbitrary number. The flux globalization based WB PCCU schemes are obtained by applying the CU numerical flux to
the semi-discretization of \eref{1.6} and by using the path-conservative technique to the evaluation of the integrals in \eref{1.7}. These
schemes have been applied to a variety of hyperbolic systems in \cite{CKL23,CKLX_22,CKLZ_23,CKN_23,KLX_21}. However, all of these WB PCCU
schemes are only second-order accurate, which results in a limited resolution of certain complicated solution structures unless very fine
mesh is used. The goal of this paper is to extend the flux globalization based WB PCCU schemes to the fifth order of accuracy via an
alternative weighted essentially non-oscillatory (A-WENO) framework.

In the past few decades, high-order finite-volume (FV) and finite-difference (FD) weighted essentially non-oscillatory (WENO) schemes have
gained their popularity due to their ability to achieve high resolution. These schemes are based on high-order WENO reconstructions and
interpolations; see the review papers \cite{Shu09,Shu_Acta} and referenced therein. A drawback of FV WENO schemes is that they use FV WENO
reconstructions, which are quite complicated and computationally expensive, especially in the multidimensional case. On the other hand, FD
WENO schemes use 1-D WENO interpolations even in the multidimensional case, which makes them substantially simpler than their FV
counterparts. The simplest and yet very accurate and powerful FD WENO schemes are A-WENO schemes were proposed in \cite{JSZ}; see also
\cite{CCK23_Adaptive,Liu17,WDGK_20,WLGD_18,WDKL}. A-WENO schemes employ standard FV numerical fluxes (with no flux splitting required) and
achieve high order of accuracy thanks to the higher-order correction terms. A-WENO schemes based on the CU numerical fluxes have been
introduced in \cite{WDGK_20,WDKL} and extended to nonconservative systems in \cite{Chu21,CKMZ23}.

In this paper, we develop flux globalization based WB A-WENO schemes that employ CU numerical fluxes and fifth-order path-conservative
integration techniques. Our A-WENO schemes utilize new higher-order correction terms, which we have recently introduced in \cite{CKX23}. The
designed schemes are then applied to the nozzle flow system and the two-layer shallow water equations. We perform a series of numerical
experiments, which clearly demonstrate that the proposed fifth-order scheme achieves higher resolution than its second-order counterpart.

The paper is organized as follows. In \S\ref{sec2}, we give a brief overview of flux globalization based WB PCCU schemes. In \S\ref{sec3},
we construct new flux globalization based WB A-WENO schemes and then apply them to the nozzle flow system in \S\ref{sec31} and the
two-layer shallow water equations in \S\ref{sec32}. In \S\ref{sec4}, we test the new schemes on a number of challenging numerical examples.
Finally, in \S\ref{sec5}, we give some concluding remarks.

\section{Flux Globalization Based WB PCCU Schemes: an Overview}\label{sec2}
In this section, we give an overview of second-order flux globalization based WB PCCU schemes, which were introduced in \cite{KLX_21} for
general nonconservative systems \eref{1.1} written in an equivalent quasi-conservative form \eref{1.6}--\eref{1.7}.

We first introduce a uniform mesh consisting of the FV cells $C_j:=[x_\jmh,x_\jph]$ of size $x_\jph-x_\jmh\equiv\dx$ centered at
$x_j=(x_\jmh+x_\jph)/2$, $j=1,\ldots,N$. We assume that at a certain time level $t$, the approximate solution, realized in terms of the cell
averages
\begin{equation*}
\xbar{\bm U}_j(t):\approx\frac{1}{\dx}\int\limits_{C_j}\bm U(x,t)\,{\rm d}x,
\end{equation*}
is available (in the rest of the paper, we will suppress the time-dependence of all of the indexed quantities for the sake of brevity). We
then evolve the cell averages $\,\xbar{\bm U}_j$ in time by solving the following system of ODEs:
\begin{equation}
\frac{{\rm d}}{{\rm d}t}\,\xbar{\bm U}_j=-\frac{\bm{{\cal K}}^{\rm FV}_\jph-\bm{{\cal K}}^{\rm FV}_\jmh}{\dx},
\label{2.1.1}
\end{equation}
where $\bm{{\cal K}}^{\rm FV}_\jph$ are the WB PCCU numerical fluxes given by (see \cite{KLX_21})
\begin{equation}\resizebox{.9\hsize}{!}{$
\bm{{\cal K}}^{\rm FV}_\jph=\frac{a^+_\jph\bm K(\bm U^-_\jph)-a^-_\jph\bm K(\bm U^+_\jph)}{a^+_\jph-a^-_\jph}+
\frac{a^+_\jph a^-_\jph}{a^+_\jph-a^-_\jph}\left(\widehat{\bm U}^+_\jph-\widehat{\bm U}^-_\jph\right),
\label{2.1.2}$}
\end{equation}
$\bm U^\pm_\jph$ and $\widehat{\bm U}^\pm_\jph$ are two slightly different approximations of the one-sided point values of
$\bm U(x_\jph,t)$, and $a^\pm_\jph$ are the one-sided local speeds of propagation, which can be estimated using the eigenvalues
$\lambda_1<\cdots<\lambda_d$ of the matrix ${\cal A}(\bm U):=\nicefrac{\partial\bm F}{\partial\bm U}-B(\bm U)$:
\begin{equation*}
\begin{aligned}
& a^-_\jph=\min\left\{\lambda_1\big({\cal A}(\bm U^-_\jph)\big),\lambda_1\big({\cal A}(\bm U^+_\jph\big)),0\right\},\\
& a^+_\jph=\max\left\{\lambda_d\big({\cal A}(\bm U^-_\jph)\big),\lambda_d\big({\cal A}(\bm U^+_\jph)\big),0\right\}.
\end{aligned}
\end{equation*}

Since $\bm K({\bm U})$ is a global flux, one needs to use a certain quadrature to evaluate its global part $\bm R(\bm U)$. The quadrature
should rely on the path-conservative technique and it also must be WB. To this end, we follow the approach introduced in \cite{KLX_21} and
compute
\begin{equation}
\bm K(\bm U^\pm_\jph)=\bm F(\bm U^\pm_\jph)-\bm R(\bm U^\pm_\jph)
\label{2.1.4}
\end{equation}
as follows.

First, in order to ensure the WB property of the resulting scheme, the point values $\bm U^\pm_\jph$ are to be obtained with the help of a
piecewise linear reconstruction of the equilibrium variables $\bm E$:
\begin{equation}
\bm E^-_\jph=\bm E_j+\frac{\dx}{2}(\bm E_x)_j,\quad\bm E^+_\jph=\bm E_{j+1}-\frac{\dx}{2}(\bm E_x)_{j+1},
\label{2.1.4a}
\end{equation}
where $\bm E_j:=\bm E(\,\xbar{\bm U}_j)$ and the slopes $(\bm E_x)_j$ are computed using a nonlinear limiter, for instance, the generalized
minmod one (see \cite{LN,NT,Swe}):
\begin{equation*}
(\bm E_x)_j={\rm minmod}\left(\theta\,\frac{\bm E_{j+1}-\bm E_j}{\dx},\frac{\bm E_{j+1}-\bm E_{j-1}}{2\dx},
\theta\,\frac{\bm E_j-\bm E_{j-1}}{\dx}\right),\quad\theta\in[1,2],
\end{equation*}
applied in a component-wise manner. Here,
\begin{equation*}
{\rm minmod}(c_1,c_2,\ldots)=
\begin{cases}
  \min(c_1,c_2,\ldots), & \mbox{if }~c_i>0,~\forall i, \\
  \max(c_1,c_2,\ldots), & \mbox{if }~c_i<0,~\forall i, \\
  0, & \mbox{otherwise}.
\end{cases}
\end{equation*}
and the parameter $\theta$ is used to control the amount of numerical dissipation present in the resulting scheme: larger values of $\theta$
correspond to sharper but, in general, more oscillatory reconstructions. The corresponding values $\bm U^\pm_\jph$ are then computed by
solving the (nonlinear systems of) equations
\begin{equation}
\bm E(\bm U^+_\jph)=\bm E^+_\jph\quad\mbox{and}\quad\bm E(\bm U^-_\jph)=\bm E^-_\jph
\label{2.1.4b}
\end{equation}
for $\bm U^+_\jph$ and $\bm U^-_\jph$, respectively.

Equipped with $\bm U^\pm_\jph$, we proceed with the evaluation of $\bm R(\bm U^\pm_\jph)$. We select $\widehat x=x_\hf$, set
$\bm R(\bm U^-_\hf):=\bm0$, evaluate
\begin{equation}
\bm R(\bm U^+_\hf)=\bm B_{\bm\Psi,\hf},
\label{2.1.5}
\end{equation}
and then recursively obtain
\begin{equation}
\bm R(\bm U^-_\jph)=\bm R(\bm U^+_\jmh)+\bm B_j,\quad\bm R(\bm U^+_\jph)=\bm R(\bm U^-_\jph)+\bm B_{\bm\Psi,\jph},
\label{2.1.6}
\end{equation}
for $j=1,\ldots,N$. In \eref{2.1.5} and \eref{2.1.6}, $\bm B_{\bm\Psi,\jph}$ and $\bm B_j$ are obtained using a proper quadrature for the integrals in
\begin{equation*}
\bm B_{\bm\Psi,\jph}=\int\limits^1_0B\big(\bm\Psi_\jph(s)\big)\bm\Psi'_\jph(s)\,{\rm d}s,\quad
\bm B_j=\int\limits_{C_j}B(\bm U)\bm U_x\,{\rm d}x,
\end{equation*}
where $\bm\Psi_\jph(s)=\bm\Psi_\jph(s;\bm U^-_\jph,\bm U^+_\jph)$ is a sufficiently smooth path connecting the states $\bm U^-_\jph$ and
$\bm U^+_\jph$, that is,
\begin{equation*}
\begin{aligned}
&\bm\Psi_\jph:[0,1]\times\mathbb R^d\times\mathbb R^d\to\mathbb R^d,\\
&\bm\Psi_\jph(0;\bm U^-_\jph,\bm U^+_\jph)=\bm U^-_\jph,\quad
\bm\Psi_\jph(1;\bm U^-_\jph,\bm U^+_\jph)=\bm U^+_\jph.
\end{aligned}
\end{equation*}
In order to ensure the resulting scheme to be WB, one has to connect the left and right equilibrium states $\bm E^-_\jph$ and
$\bm E^+_\jph$. We use a linear segment path
\begin{equation*}
\bm E_\jph(s)=\bm E^-_\jph+s\big(\bm E^+_\jph-\bm E^-_\jph\big),\quad s\in[0,1],
\end{equation*}
which is then used together with \eref{1.4} and trapezoidal quadrature to end up with
\begin{equation*}
\bm B_{\bm\Psi,\jph}=\bm F_\jph^+-\bm F_\jph^--\hf\Big[M\big(\bm U_\jph^+\big)+M\big(\bm U_\jph^-\big)\Big]
\big(\bm E_\jph^+-\bm E_\jph^-\big).
\end{equation*}
For $\bm B_j$, we similarly have
\begin{equation}
\begin{aligned}
\bm B_j&=\bm F(\bm U^-_\jph)-\bm F(\bm U^+_\jmh)-\int\limits_{C_j}M(\bm U)\bm E(\bm U)_x\,{\rm d}x\\
&\hspace{-0.6cm} \approx\bm F(\bm U^-_\jph)-\bm F(\bm U^+_\jmh)-\hf\Big[M\big(\bm U_\jph^-\big)+M\big(\bm U_\jmh^+\big)\Big]
\big(\bm E_\jph^--\bm E_\jmh^+\big).
\label{2.1.7a}
\end{aligned}
\end{equation}

Finally, in \eref{2.1.2}, $\widehat{\bm U}^\pm_\jph$ are obtained by solving modified versions of \eref{2.1.4b}:
\begin{equation*}
\widehat{\bm E}(\widehat{\bm U}^+_\jph)=\bm E^+_\jph\quad\mbox{and}\quad\widehat{\bm E}(\widehat{\bm U}^-_\jph)=\bm E^-_\jph.
\end{equation*}
Here, the functions $\widehat{\bm E}$ are the modifications of the functions $\bm E$, which are made in such a way that
$\widehat{\bm U}^+_\jph=\widehat{\bm U}^-_\jph$ as long as $\widehat{\bm E}^+_\jph=\widehat{\bm E}^-_\jph$ and, at the same time,
$\widehat{\bm U}^+_\jph=\bm U^+_\jph+{\cal O}((\dx)^2)$ and $\widehat{\bm U}^-_\jph=\bm U^-_\jph+{\cal O}((\dx)^2)$ for smooth solutions.

\section{Flux Globalization Based WB A-WENO Schemes}\label{sec3}
In this section, we extend the second-order WB PCCU schemes from \S\ref{sec2} to the fifth order of accuracy via the A-WENO framework:
\begin{equation}
\frac{{\rm d}\bm U_j}{{\rm d}t}=-\frac{\bm{{\cal K}}_\jph-\bm{{\cal K}}_\jmh}{\dx},
\label{3.1}
\end{equation}
where $\bm U_j\approx\bm U(x_j,t)$ and $\bm{{\cal K}}_\jph$ are the fifth-order A-WENO numerical fluxes:
\begin{equation*}
\bm{{\cal K}}_\jph=\bm{{\cal K}}^{\rm FV}_\jph-\frac{(\dx)^2}{24}(\bm K_{xx})_\jph+\frac{7(\dx)^4}{5760}(\bm K_{xxxx})_\jph.
\end{equation*}
Here, $\bm{{\cal K}}^{\rm FV}_\jph$ is the WB PCCU numerical flux introduced in \S\ref{sec2}, and the higher-order correction terms
$(\bm K_{xx})_\jph$ and $({\bm K_{xxxx}})_\jph$ are FD approximations of the second- and fourth-order spatial derivatives of $\bm K$ at
$x=x_\jph$, respectively. We follow \cite{CKX23} and compute the correction terms using the numerical fluxes, which have been already
obtained:
\begin{equation*}
\begin{aligned}
&(\mK_{xx})_\jph=\frac{1}{12(\dx)^2}\Big[-\bm{{\cal K}}^{\rm FV}_{j-\frac{3}{2}}+16\bm{{\cal K}}^{\rm FV}_\jmh-
30\bm{{\cal K}}^{\rm FV}_\jph+16\bm{{\cal K}}^{\rm FV}_{j+\frac{3}{2}}-\bm{{\cal K}}^{\rm FV}_{j+\frac{5}{2}}\Big],\\
&(\mK_{xxxx})_\jph=\frac{1}{(\dx)^4}\Big[\bm{{\cal K}}^{\rm FV}_{j-\frac{3}{2}}-4\bm{{\cal K}}^{\rm FV}_\jmh+6\bm{{\cal K}}^{\rm FV}_\jph-
4\bm{{\cal K}}^{\rm FV}_{j+\frac{3}{2}}+\bm{{\cal K}}^{\rm FV}_{j+\frac{5}{2}}\Big].
\end{aligned}
\end{equation*}
Finally, in order the resulting scheme to be fifth-order accurate, we apply the fifth-order affine-invariant WENO-Z (Ai-WENO-Z)
interpolation \cite{DLWW22,WD22,LLWDW23} to evaluate $\bm E^\pm_\jph$ (instead of the piecewise linear reconstruction \eref{2.1.4a}). In
addition, the integral in \eref{2.1.7a} needs to be evaluated using a quadrature, which is at least fifth-order accurate. In this paper, we
use the Newton-Cotes quadrature introduced in \cite[(4.4)]{Chu21}.

\subsection{Application to the Nozzle Flow System}\label{sec31}
In this section, we apply the flux globalization based WB A-WENO scheme to the nozzle flow system, which reads as \eref{1.1} with
\begin{equation}
\bm U=\begin{pmatrix}\sigma\rho\\\sigma\rho u\\\sigma\end{pmatrix},\quad\bm F(\bm U)=\begin{pmatrix}\sigma\rho u\\\sigma\rho u^2+\sigma p\\
0\end{pmatrix},\quad B(\bm U)=\begin{pmatrix}0&0&0\\0&0&p\\0&0&0\end{pmatrix},
\label{3.1.1}
\end{equation}
where $p(\rho)=\kappa\rho^\gamma$,
and admits steady-state solutions satisfying \eref{1.4} with
\begin{equation*}
\begin{aligned}
&M(\bm U):=\begin{pmatrix}1&0&0\\u&\sigma\rho&0\\0&0&1\end{pmatrix},\quad\bm E(\bm U)=\begin{pmatrix}q\\E\\0\end{pmatrix},\\[2.ex]
&q:=\sigma\rho u,\quad E:=\frac{u^2}{2}+\frac{\kappa\gamma}{\gamma-1}\rho^{\gamma-1}.
\end{aligned}
\end{equation*}
Here, $\rho$ is the density, $u$ is the velocity, $p$ is the pressure, $\kappa>0$ and $1<\gamma<\frac{5}{3}$ are constants, and
$\sigma=\sigma(x)$ denotes the cross-section of the nozzle, which is a given function independent of time. The nozzle flow system
\eref{1.1}, \eref{3.1.1} can be rewritten in the quasi-conservative form \eref{1.6}--\eref{1.7} with
\begin{equation}
\bm K(\bm U)=\begin{pmatrix}\sigma\rho u\\K\\0\end{pmatrix},\quad\bm R(\bm U)=\begin{pmatrix}0\\R\\0\end{pmatrix},
\label{3.1.2}
\end{equation}
where 
$$
K=\sigma\rho u^2+\sigma p-R,\quad R:=\int\limits^x_{\widehat x}p(\xi,t)\sigma_\xi(\xi)\,{\rm d}\xi.
$$
In order to construct the fifth-order WB A-WENO scheme, we first compute
\begin{equation*}
E_j=\frac{u_j^2}{2}+\frac{\kappa\gamma}{\gamma-1}(\rho_j)^{\gamma-1},
\end{equation*}
where $u_j:=\,\xbar q_j/\,\xbar{(\sigma\rho)}_j$, $\rho_j:=\xbar{(\sigma\rho)}_j/\sigma_j$, and $\sigma_j:=\sigma(x_j)$. Given $\,\xbar q_j$
and $E_j$, we perform the Ai-WENO-Z interpolation to obtain the point values $q^\pm_\jph$ and $E^\pm_\jph$, and then numerically solve the
nonlinear equations
\begin{equation*}
E^\pm_\jph=\frac{\big(q^\pm_\jph\big)^2}{2\big((\sigma\rho)_\jph^{\pm}\big)^2}+
\frac{\kappa\gamma}{\gamma-1}\big(\sigma_\jph^\pm\big)^{1-\gamma}\big[(\sigma\rho)^\pm_\jph\big]^{\gamma-1}
\end{equation*}
for $(\sigma\rho)^\pm_\jph$. Here, the values $\sigma^\pm_\jph$ are obtained by a Ai-WENO-Z interpolation of $\sigma$. For solution details,
see \cite{KLX_21}.

After obtaining the values $(\sigma\rho)^\pm_\jph$, we compute the corresponding values
$u^\pm_\jph=\nicefrac{q^\pm_\jph}{(\sigma\rho)^\pm_\jph}$ and $\rho^\pm_\jph=\nicefrac{(\sigma\rho)^\pm_\jph}{\sigma^\pm_\jph}$, which are
then used to estimate the one-sided local speeds of propagation:
\begin{equation*}
\begin{aligned}
&a^+_\jph=\max\left\{u^+_\jph+\sqrt{\kappa\gamma}\,\big(\rho^+_\jph\big)^{\frac{\gamma-1}{2}},\,
u^-_\jph+\sqrt{\kappa\gamma}\,\big(\rho^-_\jph\big)^{\frac{\gamma-1}{2}},\,0\right\},\\
&a^-_\jph=\min\left\{u^+_\jph-\sqrt{\kappa\gamma}\,\big(\rho^+_\jph\big)^{\frac{\gamma-1}{2}},\,
u^-_\jph-\sqrt{\kappa\gamma}\,\big(\rho^-_\jph\big)^{\frac{\gamma-1}{2}},\,0\right\}.
\end{aligned}
\end{equation*}
The global flux $\bm{K}$ in \eref{3.1.2} needs to be evaluated by \eref{2.1.4}, which for the nozzle flow system reads as
\begin{equation*}
\bm K(\bm U^\pm_\jph)=\bm F(\bm U^\pm_\jph)-\bm R(\bm U^\pm_\jph)=
\begin{pmatrix}q^\pm_\jph\\q^\pm_\jph u^\pm_\jph+\sigma^\pm_\jph p^\pm_\jph\\0\end{pmatrix}-\begin{pmatrix}0\\R^\pm_\jph\\0\end{pmatrix},
\end{equation*}
where $p^\pm_\jph=\kappa\big(\rho^\pm_\jph\big)^\gamma$ and $R^\pm_\jph$ are computed as in \S\ref{sec2}. In particular, we have
$R^-_\hf:=0$, $R^+_\hf=B_{\bm\Psi,\hf}$, and
\begin{equation*}
R^-_\jph=R^+_\jmh+B_j,\quad R^+_\jph=R^-_\jph+B_{\bm\Psi,\jph},\quad j=1,\ldots,N,
\end{equation*}
where
\begin{equation*}
\begin{aligned}
B_{\bm\Psi,\jph}&=q^+_\jph u^+_\jph+\sigma^+_\jph p^+_\jph-q^-_\jph u^-_\jph-\sigma^-_\jph p^-_\jph\\
&-\frac{u^+_\jph+u^-_\jph}{2}\big(q^+_\jph-q^-_\jph\big)-\frac{(\sigma\rho)^+_\jph+(\sigma\rho)^-_\jph}{2}\big(E^+_\jph-E^-_\jph\big),
\end{aligned}
\end{equation*}
and
\begin{equation*}
B_j=q^-_\jph u^-_\jph+\sigma^-_\jph p^-_\jph-q^+_\jmh u^+_\jmh-\sigma^+_\jmh p^+_\jmh-\int\limits_{C_j}(uq_x+\sigma\rho E_x)\,{\rm d}x,
\end{equation*}
where the integral in the last formula is evaluated using the fifth-order Newton-Cotes quadrature \cite[(4.4)]{Chu21}.

Finally, $\widehat{\bm U}^\pm_\jph:=((\widehat{\sigma\rho})^{\,\pm}_\jph,\widehat q^{\,\pm}_\jph)^\top$, where
$\widehat q^{\,\pm}_\jph=q^\pm_\jph$ and $(\widehat{\sigma\rho})^\pm_\jph$ are obtained by numerically solving the following nonlinear
equations:
\begin{equation*}
E^\pm_\jph=\frac{\big(q^\pm_\jph\big)^2}{2\big((\widehat{\sigma\rho})_\jph^\pm\big)^2}+
\frac{\kappa\gamma}{\gamma-1}\big(\widehat\sigma_\jph\big)^{1-\gamma}\big[(\widehat{\sigma\rho})^\pm_\jph\big]^{\gamma-1}
\end{equation*}
with $\widehat\sigma_\jph:=\nicefrac{(\sigma_\jph^++\sigma_\jph^-)}{2}$; for solution details, see \cite{KLX_21}.

\subsection{Application to the Two-Layer Shallow Water Equations}\label{sec32}
In this section, we apply the flux globalization based WB A-WENO scheme to the two-layer shallow water equations, which read as \eref{1.1}
with
\begin{equation}
\begin{aligned}
\bm U&=\begin{pmatrix}h_1\\q_1\\h_2\\q_2\\Z\end{pmatrix},\quad
\bm F(\bm U)=\begin{pmatrix}q_1\\h_1u_1^2+\dfrac{g}{2}h_1^2\\q_2\\h_2u_2^2+\dfrac{g}{2}h_2^2\\0\end{pmatrix},\\
B(\bm U)&=\begin{pmatrix}0&0&0&0&0\\0&0&-gh_1&0&-gh_1\\0&0&0&0&0\\-rgh_2&0&0&0&-gh_2\\0&0&0&0&0\end{pmatrix},
\label{3.2.1}
\end{aligned}
\end{equation}
and admit steady-state solutions satisfying \eref{1.4} with
\begin{equation*}
M(\bm U):=\begin{pmatrix}1&0&0&0&0\\u_1&h_1&0&0&0\\0&0&1&0&0\\0&0&u_2&h_2&0\\0&0&0&0&1\end{pmatrix},\quad
\bm E(\bm U):=\begin{pmatrix}q_1\\E_1\\q_2\\E_2\\0\end{pmatrix},
\end{equation*}
where
\begin{equation*}
E_1:=\frac{q_1^2}{2h_1^2}+g(h_1+h_2+Z)\quad\mbox{and}\quad E_2:=\frac{q_2^2}{2h_2^2}+g(rh_1+h_2+Z).
\end{equation*}
Here, $h_1$ and $h_2$ are the water depths in the upper and lower layers, respectively, $u_1$ and $u_2$ are the corresponding velocities,
$q_1:=h_1u_1$ and $q_2:=h_2u_2$ represent the corresponding discharges, $Z(x)$ is a function describing the bottom topography, which can be
discontinuous, $g$ is the constant acceleration due to gravity, and $r:=\frac{\rho_1}{\rho_2}<1$ is the ratio of the constant densities
$\rho_1$ (upper layer) and $\rho_2$ (lower layer). The two-layer shallow water equations \eref{1.1}, \eref{3.2.1} can be rewritten in the
quasi-conservative form \eref{1.6}--\eref{1.7} with
\begin{equation*}
\bm K(\bm U)=(q_1,K_1,q_2,K_2,0)^\top,\quad\bm R(\bm U)=(0,R_1,0,R_2,0)^\top,
\end{equation*}
where
\begin{equation}
\begin{aligned}
&K_1=h_1u_1^2+\frac{g}{2}h_1^2-R_1,&&R_1:=-\int\limits^x_{\widehat x}gh_1(\xi,t)\big[h_2(\xi,t)+Z(\xi)\big]_\xi\,{\rm d}\xi,\\
&K_2=h_2u_2^2+\frac{g}{2}h_2^2-R_2,&&R_2:=-\int\limits^x_{\widehat x}gh_2(\xi,t)\big[rh_1(\xi,t)+Z(\xi)\big]_\xi\,{\rm d}\xi.
\end{aligned}
\label{3.2.2f}
\end{equation}

In order to construct the fifth-order WB A-WENO scheme, we first compute
\begin{equation*}
\begin{aligned}
(E_1)_j:&=\frac{(\xbar{q_1})_j^2}{2(\xbar{h_1})_j^2}+g\left[(\xbar{h_1})_j+(\xbar{h_2})_j+Z_j\right],\\
(E_2)_j:&=\frac{(\xbar{q_2})_j^2}{2(\xbar{h_2})_j^2}+g\left[r(\xbar{h_1})_j+(\xbar{h_2})_j+Z_j\right],
\end{aligned}
\end{equation*}
where $Z_j:=Z(x_j)$. Given $(\,\xbar{q_1})_j$, $(\,\xbar{q_2})_j$, $(E_1)_j$, and $(E_2)_j$, we perform the Ai-WENO-Z interpolation to
obtain the point values $\bm E^\pm_\jph=((q_1)^\pm_\jph,(E_1)^\pm_\jph,$ $(q_2)^\pm_\jph, (E_2)^\pm_\jph,0)^\top$, and then at every cell
interfaces $x=x_\jph$, we numerically solve the (nonlinear) systems of equations
\begin{equation}
\left\{\begin{aligned}
&(E_1)^+_\jph=\frac{\big((q_1)^+_\jph\big)^2}{2\big((h_1)^+_\jph\big)^2}+g\left[(h_1)^+_\jph+(h_2)^+_\jph+Z^+_\jph\right],\\
&(E_2)^+_\jph=\frac{\big((q_2)^+_\jph\big)^2}{2\big((h_2)^+_\jph\big)^2}+g\left[r(h_1)^+_\jph+(h_2)^+_\jph+Z^+_\jph\right],
\end{aligned}\right.
\label{3.2.4}
\end{equation}
and
\begin{equation}
\left\{\begin{aligned}
&(E_1)^-_\jph=\frac{\big((q_1)^-_\jph\big)^2}{2\big((h_1)^-_\jph\big)^2}+g\left[(h_1)^-_\jph+(h_2)^-_\jph+Z^-_\jph\right],\\
&(E_2)^-_\jph=\frac{\big((q_2)^-_\jph\big)^2}{2\big((h_2)^-_\jph\big)^2}+g\left[r(h_1)^-_\jph+(h_2)^-_\jph+Z^-_\jph\right],
\end{aligned}\right.
\label{3.2.4a}
\end{equation}
for $\big((h_1)^+_\jph,(h_2)^+_\jph\big)$ and $\big((h_1)^-_\jph,(h_2)^-_\jph\big)$, respectively. Here, the values $Z^\pm_\jph$ in
\eref{3.2.4} and \eref{3.2.4a} are obtained using the Ai-WENO-Z interpolation of $Z$. For solution details, see \cite{KLX_21}.

After obtaining the values $(h_1)^\pm_\jph$ and $(h_2)^\pm_\jph$, we compute the corresponding values
${(u_i)}^\pm_\jph=\nicefrac{{(q_i)}^\pm_\jph}{(h_i)^\pm_\jph}$, $i=1,2$, and then evaluate the global numerical fluxes
$\bm K(\bm U_\jph^\pm)$ at the cell interfaces $x=x_\jph$ using \eref{3.2.1}--\eref{3.2.2f}:
\begin{equation*}
\begin{aligned}
\bm K(\bm U_\jph^\pm)&=\bm F(\bm U_\jph^\pm)-\bm R(\bm U_\jph^\pm)\\
&=\begin{pmatrix}(q_1)^\pm_\jph\\(h_1)^\pm_\jph\big((u_1)^\pm_\jph\big)^2+\dfrac{g}{2}\big((h_1)^\pm_\jph\big)^2\\
(q_2)^\pm_\jph\\(h_2)^\pm_\jph\big((u_2)^\pm_\jph\big)^2+\dfrac{g}{2}\big((h_2)^\pm_\jph\big)^2\\0
\end{pmatrix}-\begin{pmatrix}0\\(R_1)^\pm_\jph\\0\\(R_2)^\pm_\jph\\0\end{pmatrix},
\end{aligned}
\end{equation*}
where $(R_i)^\pm_\jph$, $i=1,2$ are computed as in \S\ref{sec2}. In particular, we have $(R_i)^-_\hf:=0$, $(R_i)^+_\hf=(B_i)_{\bm\Psi,\hf}$,
and
\begin{equation*}
(R_i)^-_\jph=(R_i)^+_\jmh+(B_i)_j,\quad (R_i)^+_\jph=(R_i)^-_\jph+(B_i)_{\bm\Psi,\jph},\quad j=1,\ldots,N,
\end{equation*}
where
\begin{equation*}
\begin{aligned}
(B_i)_{\bm\Psi,\jph}&=(q_i)^+_\jph(u_i)^+_\jph-(q_i)^-_\jph(u_i)^-_\jph+
\frac{g}{2}\left[\big((h_i)^+_\jph\big)^2-\big((h_i)^-_\jph\big)^2\right]\\
&-\frac{(u_i)^+_\jph+(u_i)^-_\jph}{2}\left[(q_i)^+_\jph-(q_i)^-_\jph\right] \\
&-\frac{(h_i)^+_\jph+(h_i)^-_\jph}{2}\left[(E_i)^+_\jph-(E_i)^-_\jph\right],
\end{aligned}
\end{equation*}
and
\begin{equation*}
\begin{aligned}
(B_i)_j&=(q_i)^-_\jph(u_i)^-_\jph-(q_i)^+_\jmh(u_i)^+_\jmh+\frac{g}{2}\left[\big((h_i)^-_\jph\big)^2-\big((h_i)^+_\jmh\big)^2\right] \\
&-\int\limits_{C_j}(u_i(q_i)_x+h_i(E_i)_x)\,{\rm d}x,
\end{aligned}
\end{equation*}
$i=1,2$. As in \S\ref{sec31}, the integral in the last formula is evaluated using the fifth-order Newton-Cotes quadrature
\cite[(4.4)]{Chu21}. The one-sided local speed $a^\pm_\jph$ are estimated using the Lagrange theorem precisely as it was done in
\cite[\S2.3]{KPmultil}.

Finally, $\widehat{\bm U}^\pm_\jph:=((\widehat h_1)^\pm_\jph,(\widehat q_1)^\pm_\jph,(\widehat h_2)^\pm_\jph,(\widehat q_2)^\pm_\jph)^\top$,
where $(\widehat q_i)^\pm_\jph=(q_i)^\pm_\jph$, $i=1,2$, and $(\widehat{h_i})^\pm_\jph$, $i=1,2$, are obtained by numerically solving the
following (nonlinear) systems:
\begin{equation*}
\left\{\begin{aligned}
&(E_1)^+_\jph=\frac{((q_1)^+_\jph)^2}{2((\widehat h_1)^+_\jph)^2}+
g\left[(\widehat h_1)^+_\jph+(\widehat h_2)^+_\jph+\widehat Z_\jph\right],\\
&(E_2)^+_\jph=\frac{((q_2)^+_\jph)^2}{2((\widehat h_2)^+_\jph)^2}+
g\left[r(\widehat h_1)^+_\jph+(\widehat h_2)^+_\jph+\widehat Z_\jph\right],
\end{aligned}\right.
\end{equation*}
and
\begin{equation*}
\left\{\begin{aligned}
&(E_1)^-_\jph=\frac{((q_1)^-_\jph)^2}{2((\widehat h_1)^-_\jph)^2}+
g\left[(\widehat h_1)^-_\jph+(\widehat h_2)^-_\jph+\widehat Z_\jph\right],\\
&(E_2)^-_\jph=\frac{((q_2)^-_\jph)^2}{2((\widehat h_2)^-_\jph)^2}+
g\left[r(\widehat h_1)^-_\jph+(\widehat h_2)^-_\jph+\widehat Z_\jph\right],
\end{aligned}\right.
\end{equation*}
with $\widehat Z_\jph:=\nicefrac{(Z_\jph^++Z_\jph^-)}{2}$; for solution details, see \cite{KLX_21}.

\section{Numerical Examples}\label{sec4}
In this section, we apply the developed fifth-order WB A-WENO schemes to several numerical examples taken from \cite{KLX_21} and compare
their performance with that of the corresponding second-order schemes. For the sake of brevity, the tested schemes will be referred to as
the 5-Order Scheme and 2-Order Scheme, respectively.

In all of the examples, we have solved the ODE systems \eref{2.1.1} and \eref{3.1} using the three-stage third-order strong stability
preserving (SSP) Runge-Kutta solver (see, e.g., \cite{Gottlieb11,Gottlieb12}) and use the CFL number 0.5. The 2-Order Scheme has been
implemented with the minmod parameter $\theta=1.3$.

\subsection{Nozzle Flow System}
We begin with two examples for the nozzle flow system \eref{1.1}, \eref{3.1.1} with $\kappa=1$ and $\gamma=1.4$. In both examples, we
consider nontrivial steady states, which can be exactly preserved by both 2-Order and 5-Order Schemes at discrete level (this has been
numerically verified), and test the ability of the studied schemes to capture small perturbations of these steady states.

\subsubsection*{Example 1---Flows in Continuous Convergent and Divergent Nozzles}
In the first example, we consider convergent and divergent nozzles described using the smooth cross-sections
\begin{equation*}
\sigma(x)=0.976-0.748\tanh(0.8x-4)\quad\mbox{and}\quad\sigma(x)=0.976+0.748\tanh(0.8x-4),
\end{equation*}
respectively. 

We take the steady states with $q_{\rm eq}\equiv8$, $E_{\rm eq}=58.3367745090349$ and $q_{\rm eq}=8$, $E_{\rm eq}= 21.9230562619897$ for the
convergent and divergent nozzles, respectively, and compute the discrete values of $\rho_{\rm eq}(x)$ in both the convergent and divergent
cases by solving the corresponding nonlinear equations. We then obtain
$u_{\rm eq}(x)=\nicefrac{q_{\rm eq}(x)}{(\sigma(x)\rho_{\rm eq}(x))}$.

Equipped with these steady states, we proceed by adding a small perturbation to the density field and consider the initial data
\begin{equation*}
\rho(x,0)=\rho_{\rm eq}(x)+\begin{cases}10^{-2},&x\in[2.9,3.1],\\0,&\mbox{otherwise},\end{cases}\qquad
q(x,0)=\sigma(x)\rho(x,0)u_{\rm eq}(x),
\end{equation*}
prescribed in the computational domain $[0,10]$ subject to free boundary conditions.

\begin{figure}[ht!]
\centerline{\includegraphics[trim=0.4cm 0.3cm 0.9cm 0.1cm, clip, width=4.cm]{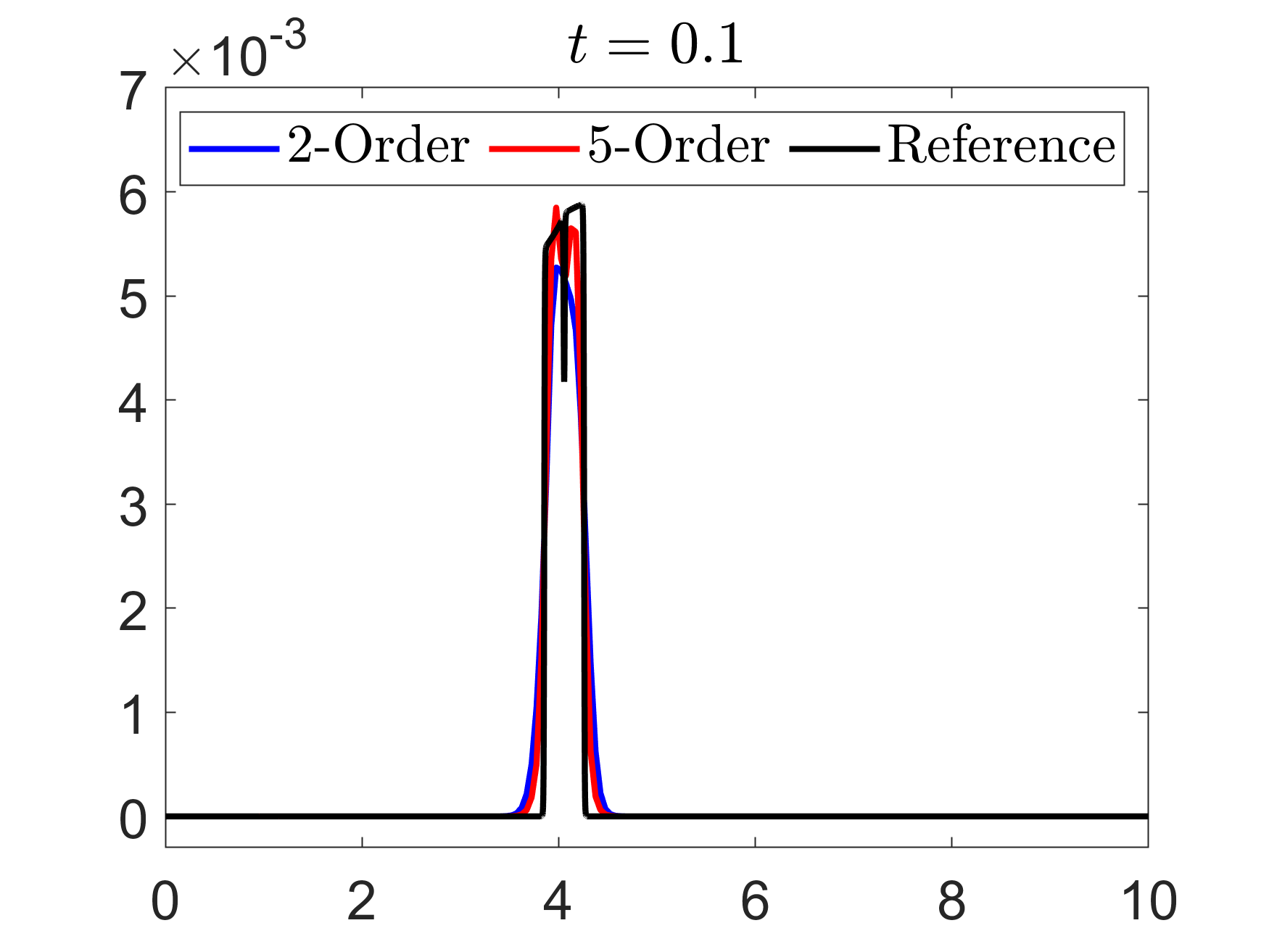}\hspace{0.1cm}
            \includegraphics[trim=0.4cm 0.3cm 0.9cm 0.1cm, clip, width=4.cm]{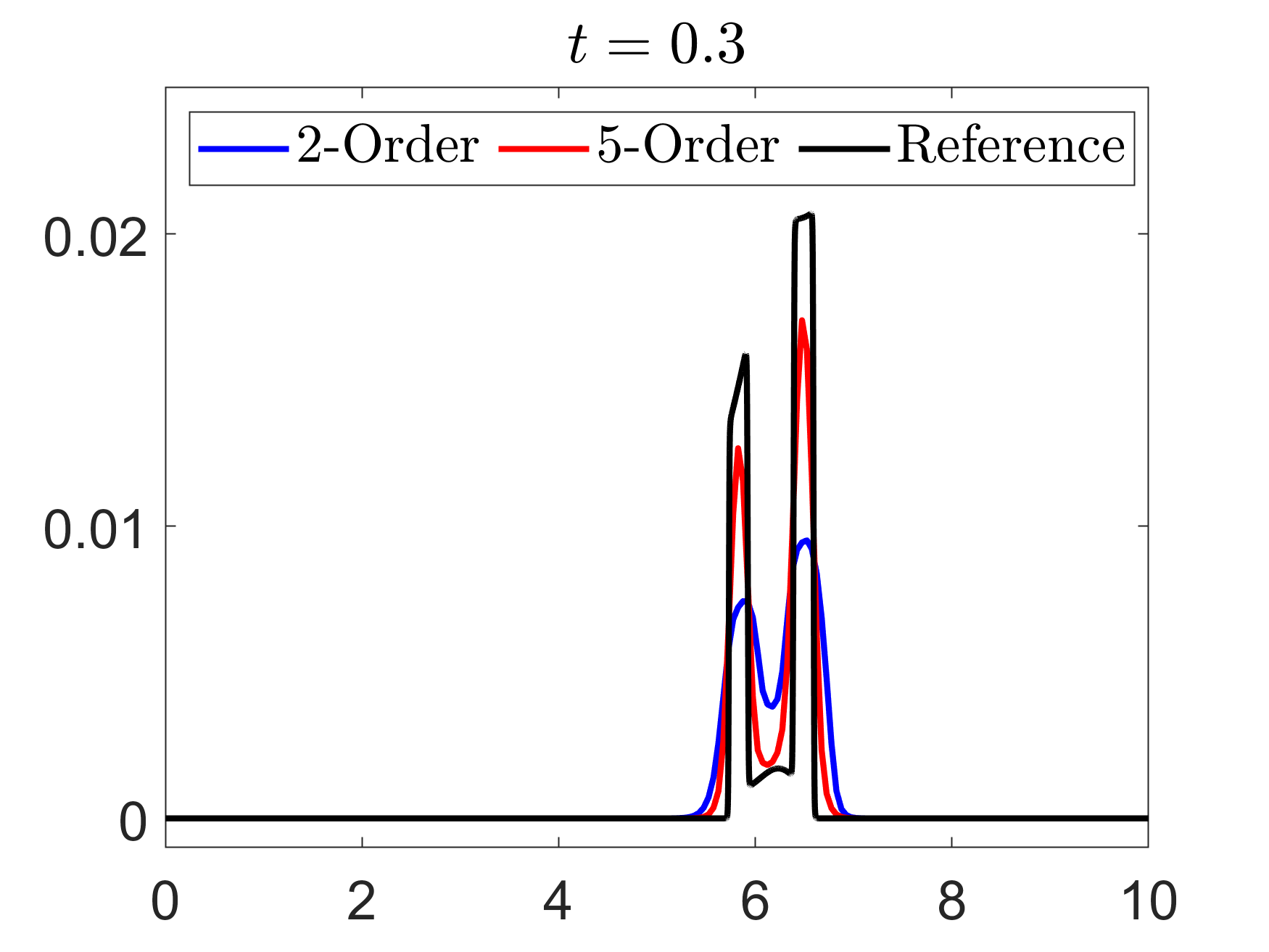}\hspace{0.1cm}
            \includegraphics[trim=0.4cm 0.3cm 0.9cm 0.1cm, clip, width=4.cm]{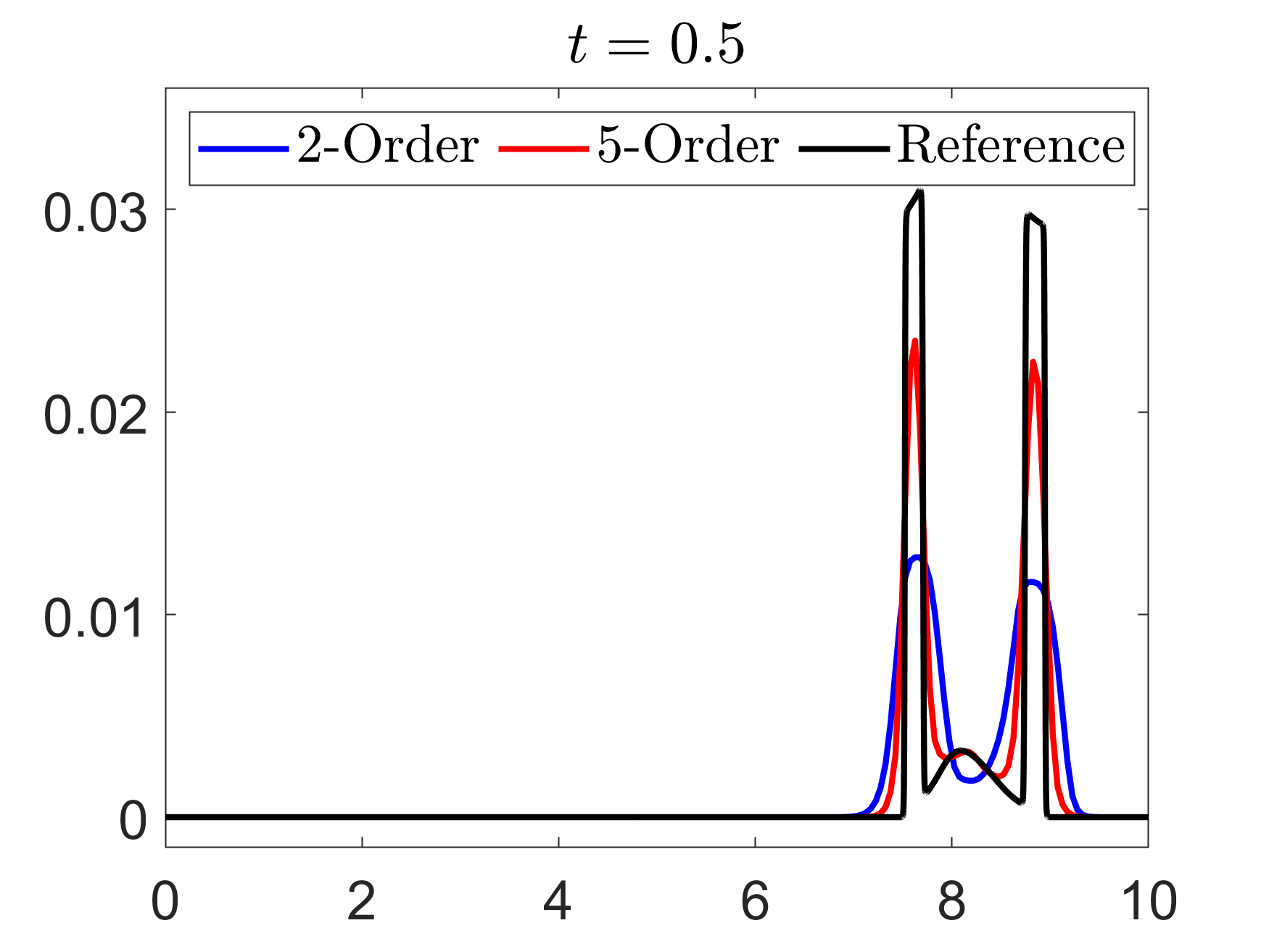}}
\vskip5pt
\centerline{\includegraphics[trim=0.4cm 0.3cm 0.9cm 0.1cm, clip, width=4.cm]{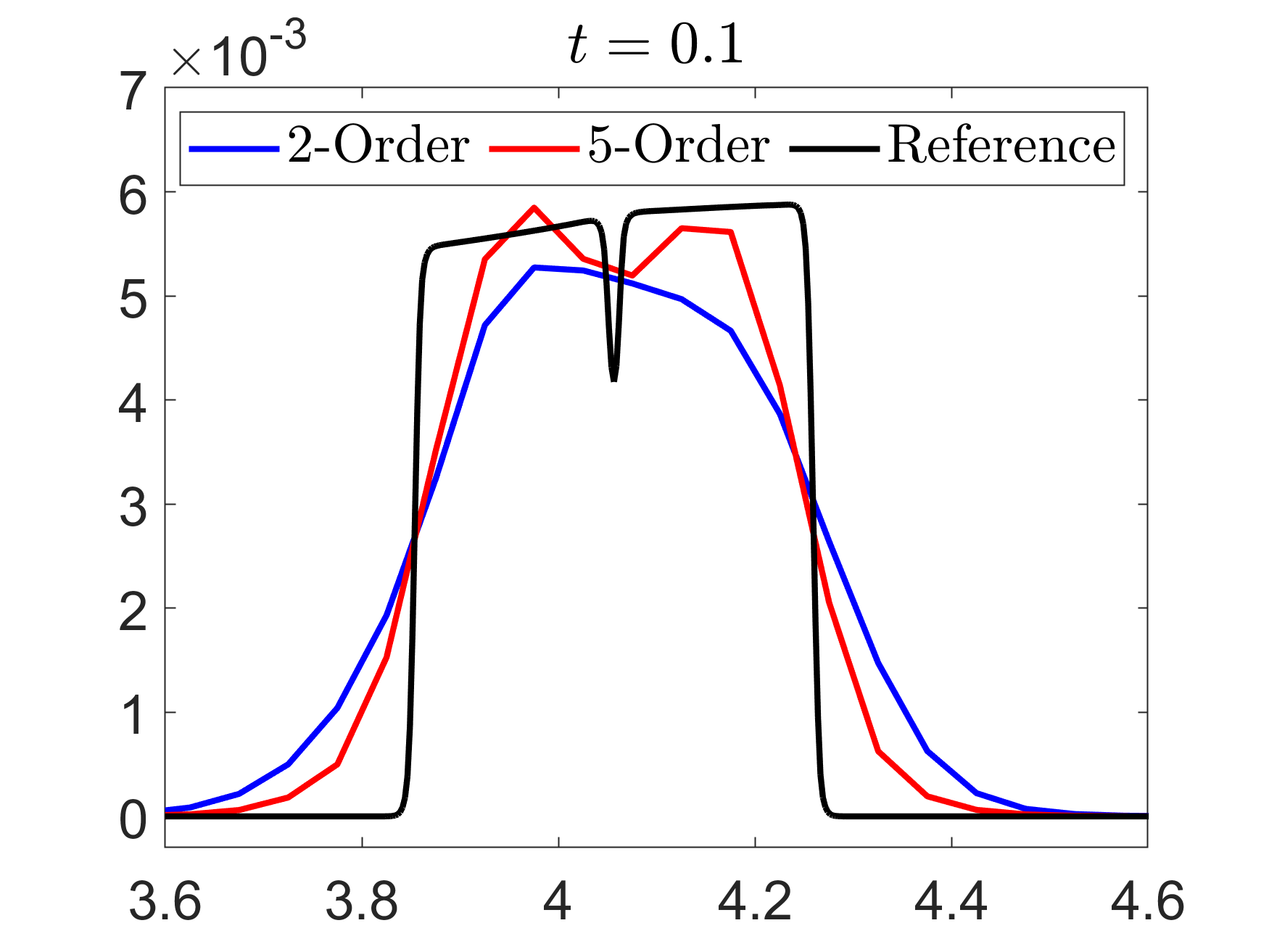}\hspace{0.1cm}
            \includegraphics[trim=0.4cm 0.3cm 0.9cm 0.1cm, clip, width=4.cm]{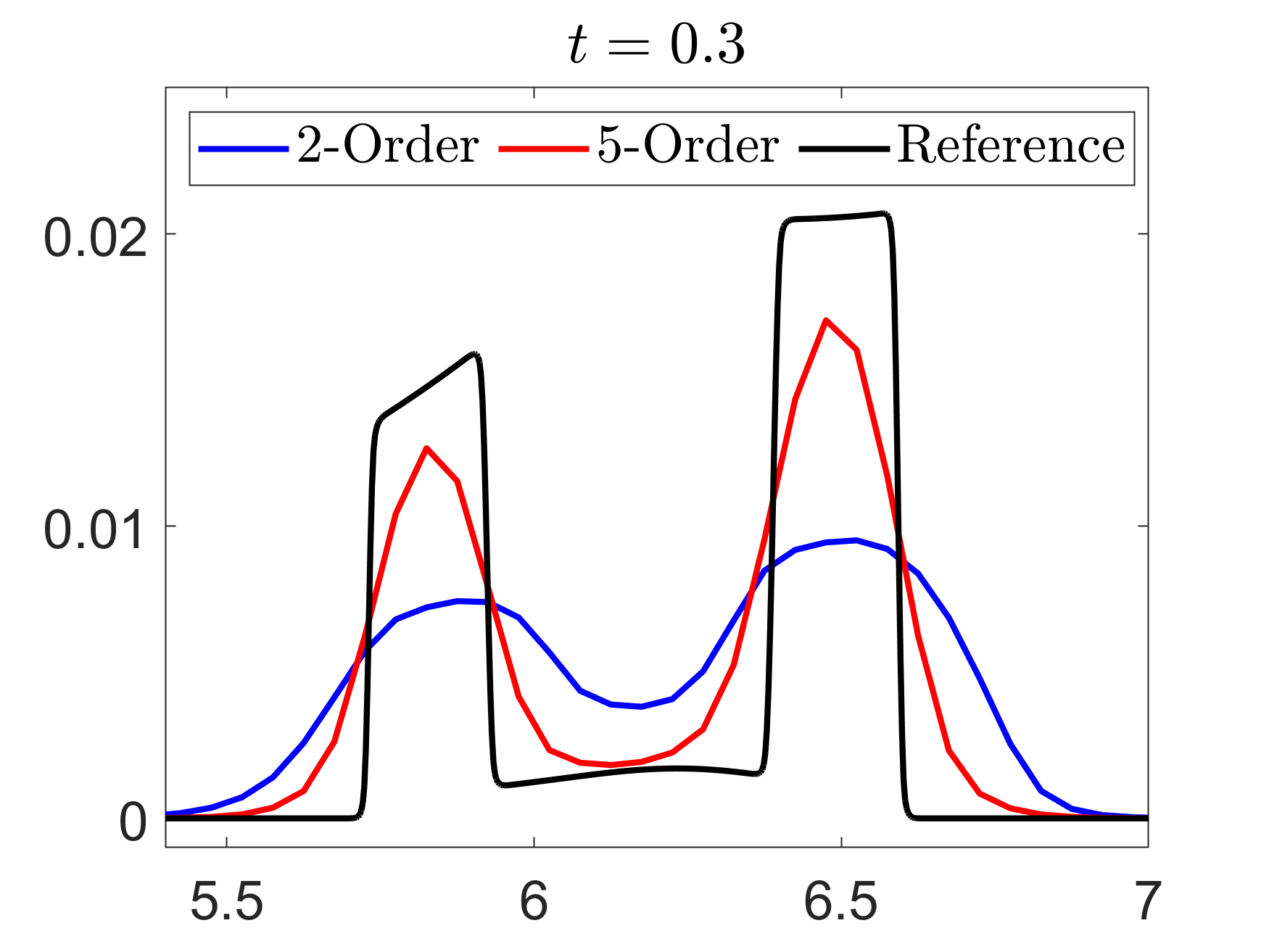}\hspace{0.1cm}
            \includegraphics[trim=0.4cm 0.3cm 0.9cm 0.1cm, clip, width=4.cm]{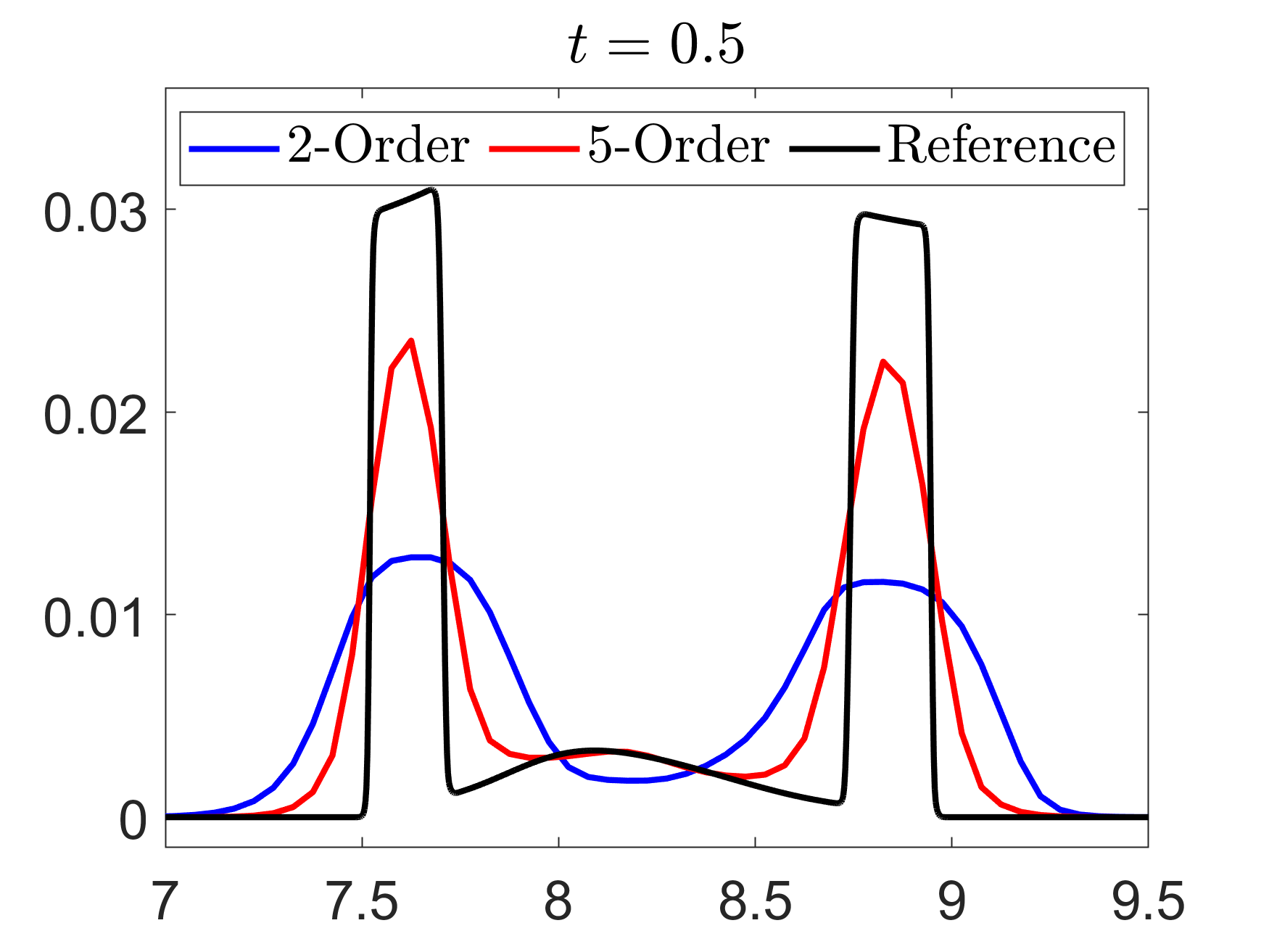}}
\caption{\sf Example 1 (Convergent Nozzle): The difference $\rho(x,t)-\rho_{eq}$ (top row) and zoom at the perturbations (bottom row) at
times $t=0.1$ (left column), 0.3 (middle column), and 0.5 (right column).\label{fig5}}
\end{figure}

\begin{figure}[ht!]
\centerline{\includegraphics[trim=0.4cm 0.3cm 0.9cm 0.1cm, clip, width=4.cm]{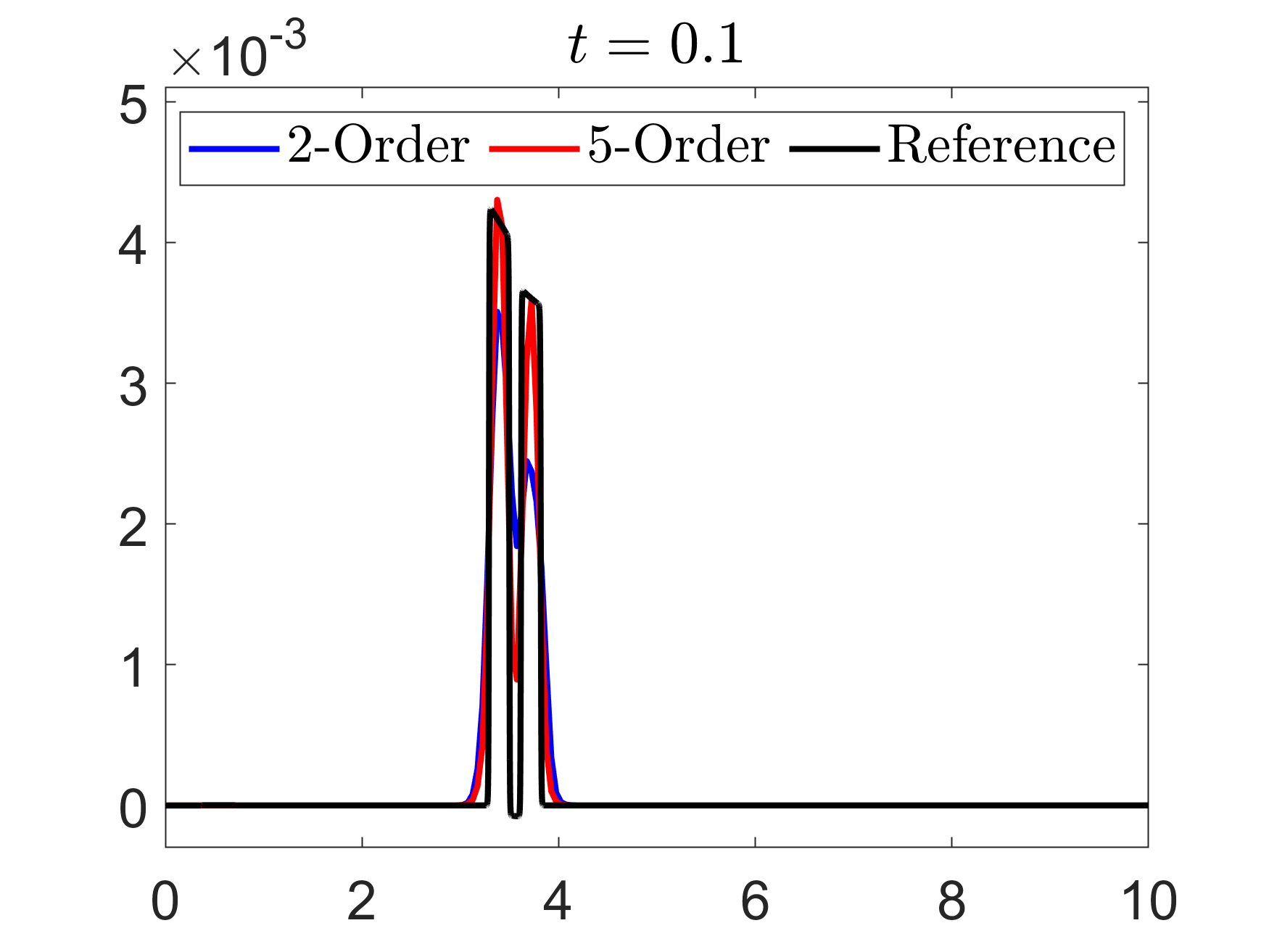}\hspace{0.1cm}
            \includegraphics[trim=0.4cm 0.3cm 0.9cm 0.1cm, clip, width=4.cm]{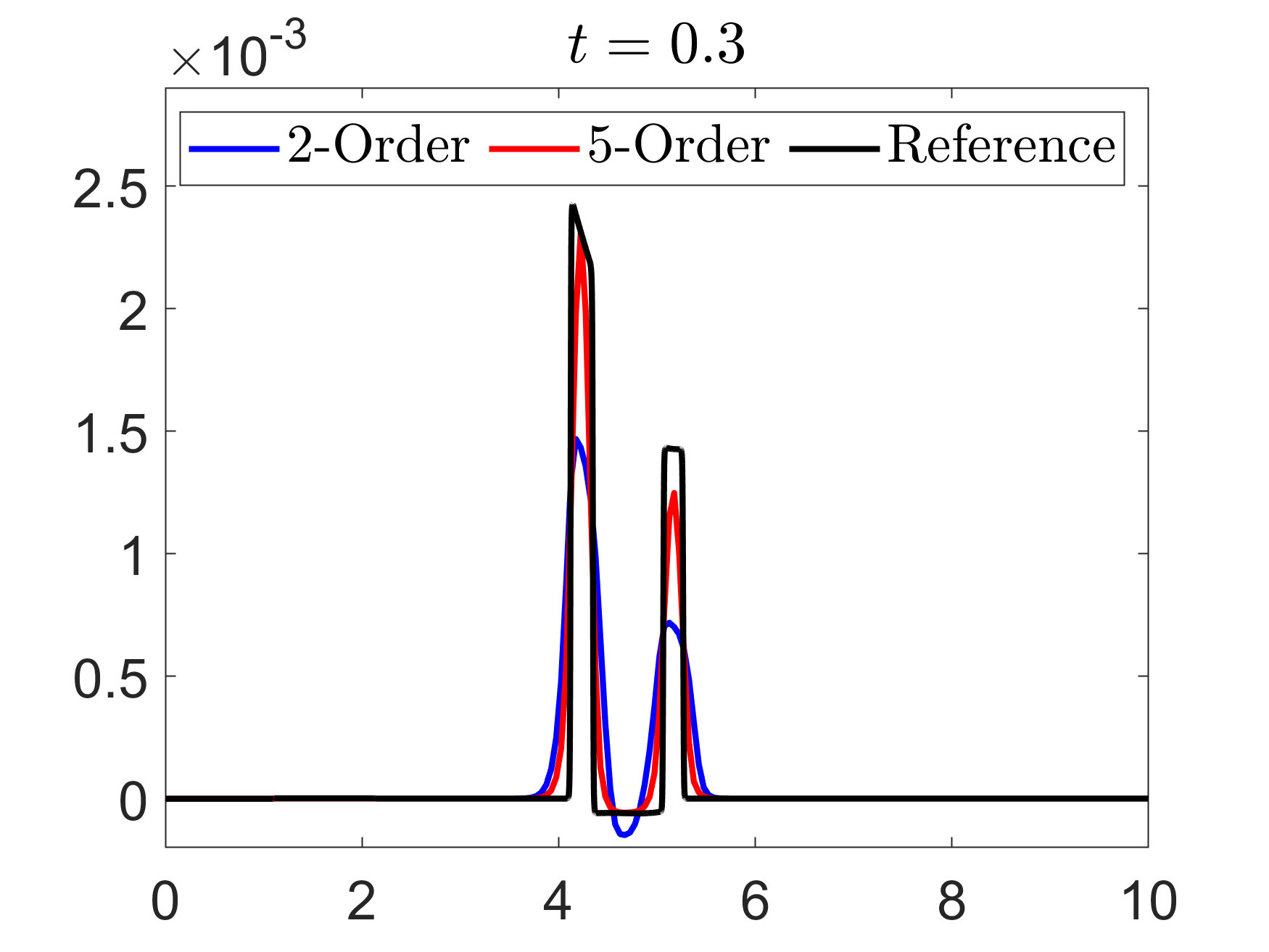}\hspace{0.1cm}
            \includegraphics[trim=0.4cm 0.3cm 0.9cm 0.1cm, clip, width=4.cm]{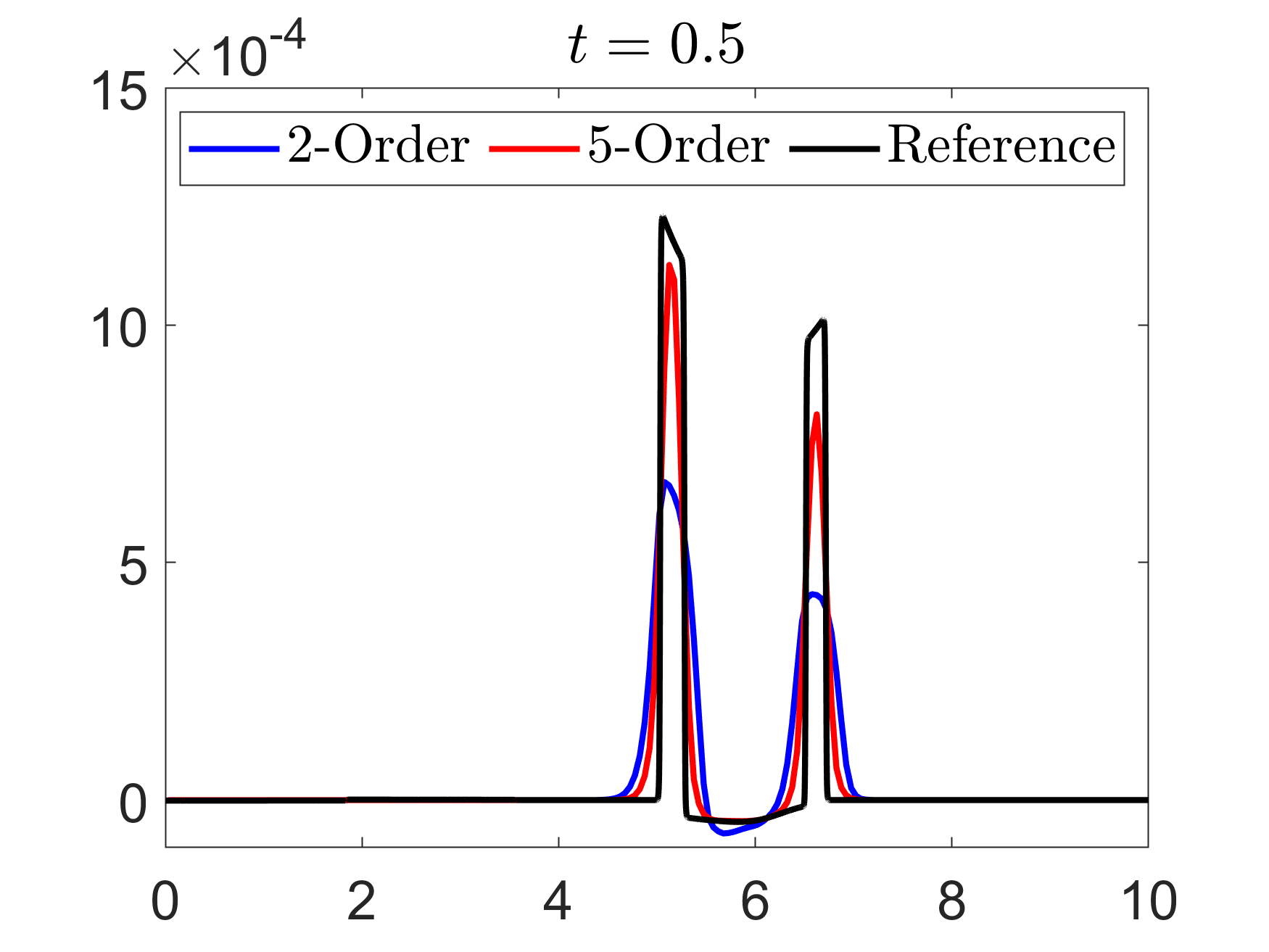}}
\vskip5pt
\centerline{\includegraphics[trim=0.4cm 0.3cm 0.9cm 0.1cm, clip, width=4.cm]{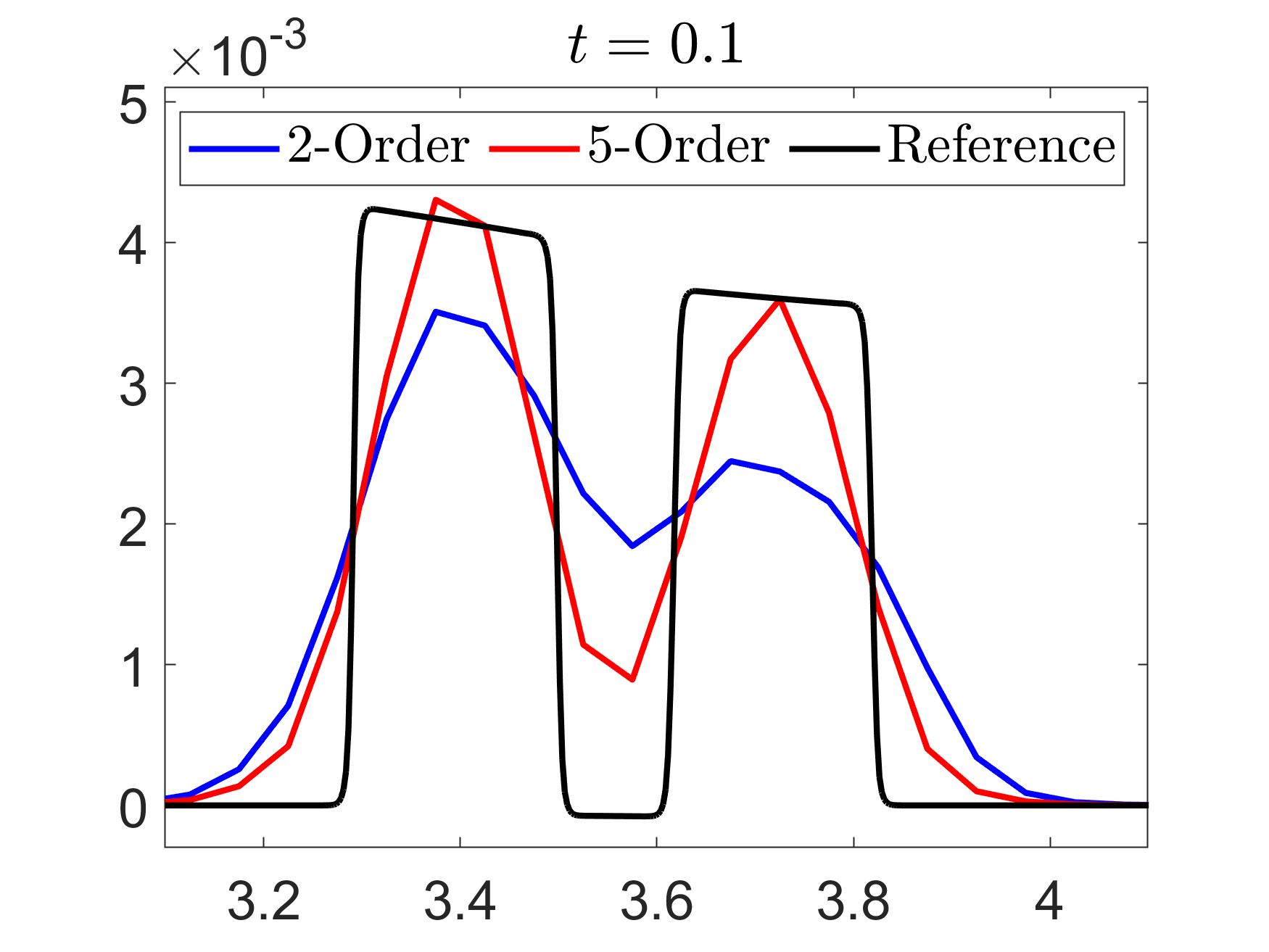}\hspace{0.1cm}
            \includegraphics[trim=0.4cm 0.3cm 0.9cm 0.1cm, clip, width=4.cm]{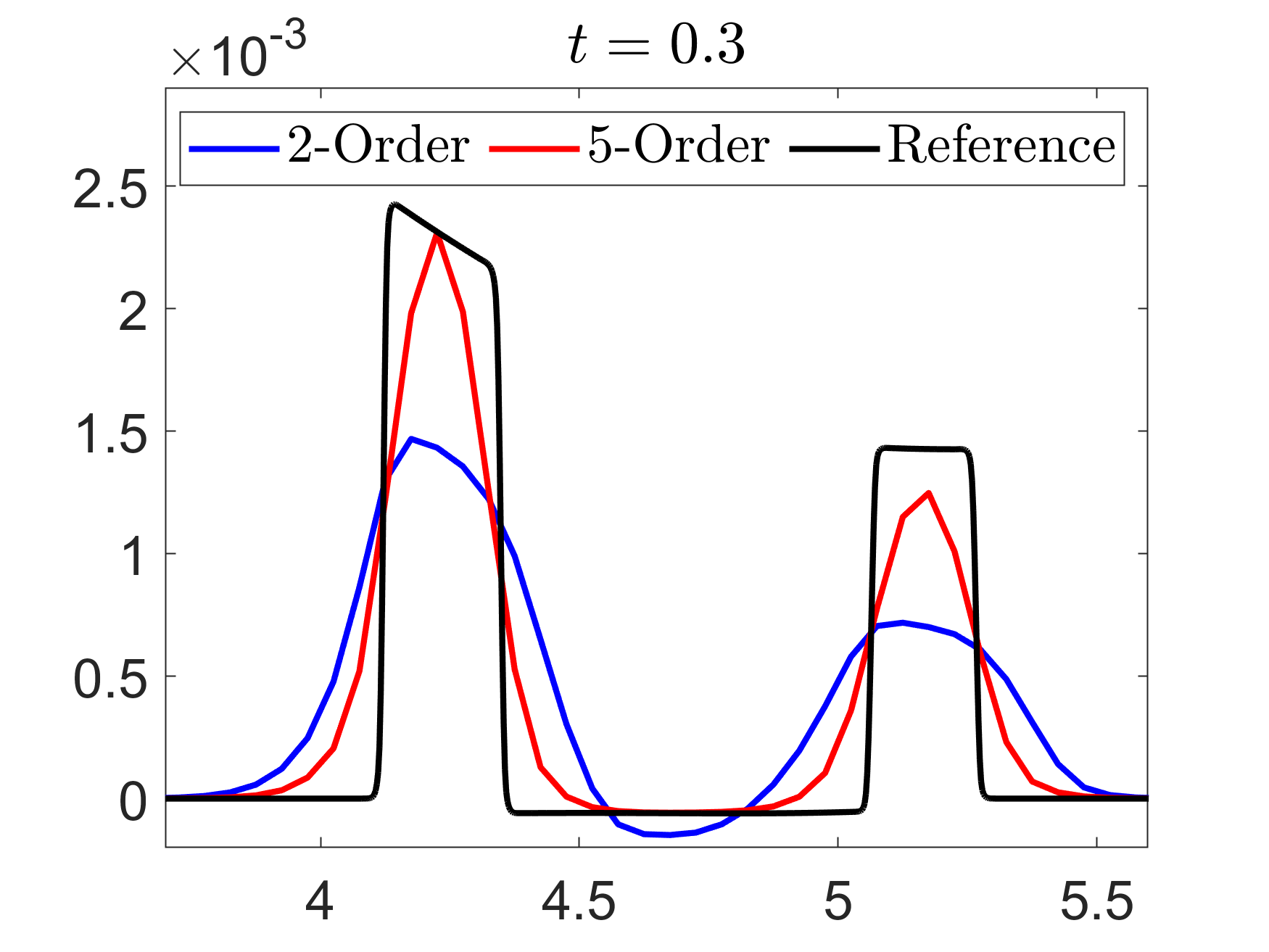}\hspace{0.1cm}
            \includegraphics[trim=0.4cm 0.3cm 0.9cm 0.1cm, clip, width=4.cm]{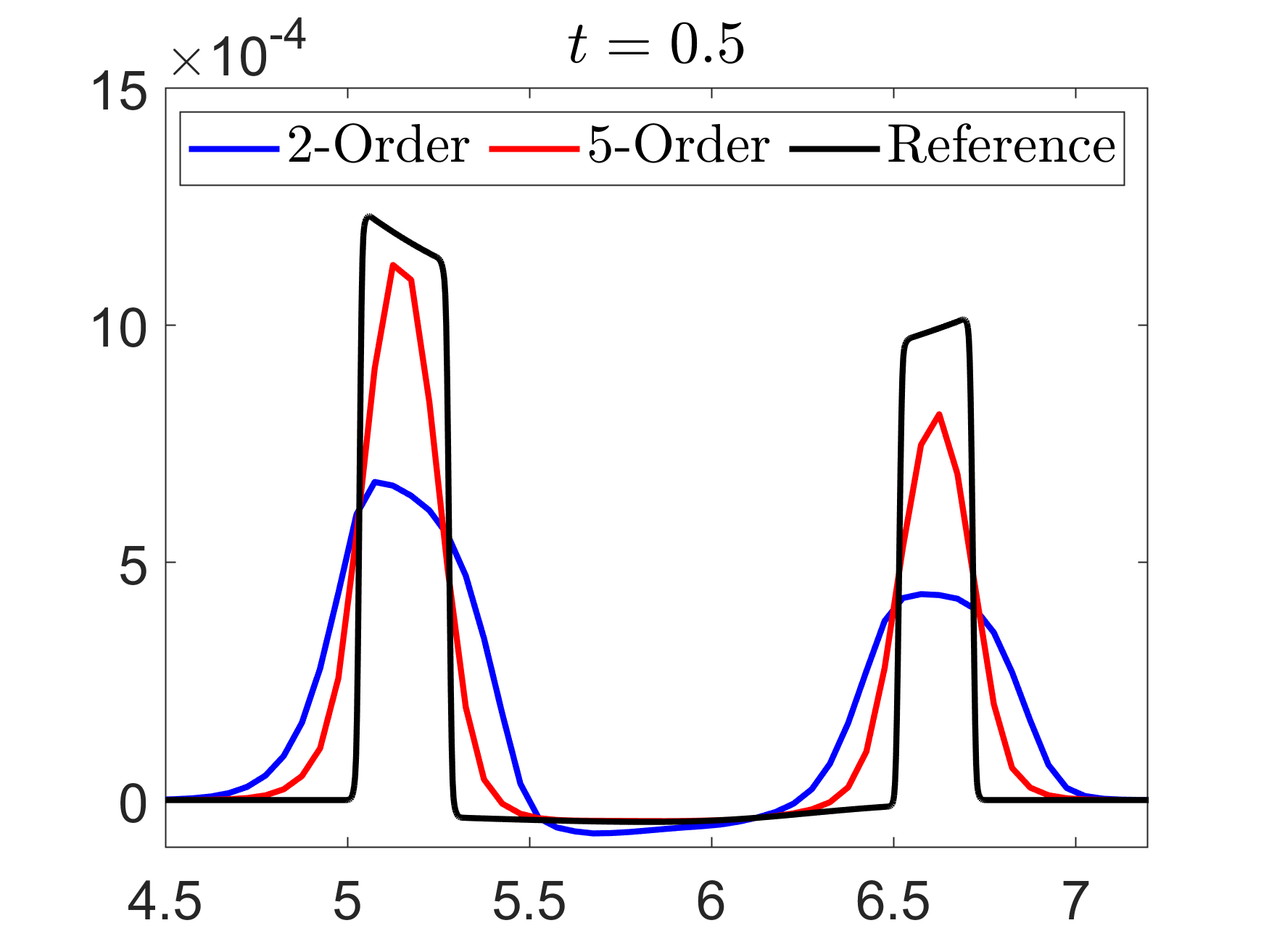}}
\caption{\sf Example 1: Same as in Figure \ref{fig5}, but for the divergent nozzle.\label{fig6}}
\end{figure}

We compute the numerical solutions until the final time $t=0.5$ by both the 2-Order and 5-Order Schemes for both convergent and divergent
nozzles on a uniform mesh with $\dx=1/20$ and plot the differences $\rho(x,t)-\rho_{\rm eq}(x)$ at times $t=0.1$, 0.3, and 0.5 in Figures
\ref{fig5} and \ref{fig6}. We also plot the reference solutions computed by the 5-Order Scheme with $\dx=1/400$. As one can
see, the 5-Order Scheme clearly outperforms the 2-Order Scheme as it achieves much higher resolution of the perturbations.

\subsubsection*{Example 2---Flows in a Discontinuous Convergent-Divergent Nozzle}
In the second example, we consider a convergent-divergent nozzle with a discontinuous $\sigma(x)$ given by
\begin{equation*}
\sigma(x)=\begin{cases}2,&x\in[0,7.5]\cup[12.5,20],\\1,&\mbox{otherwise}.\end{cases}
\end{equation*}
We take a steady state with $q_{\rm eq}\equiv8$ and $E_{\rm eq}=57.13486505$ and compute the discrete values of $\rho_{\rm eq}(x)$ by
solving the corresponding nonlinear equations. We then obtain $u_{\rm eq}(x)=\nicefrac{q_{\rm eq}(x)}{(\sigma(x)\rho_{\rm eq}(x))}$.

Equipped with these steady states, we proceed by adding a small perturbation to the density field and consider the initial data
\begin{equation*}
\rho(x,0)=\rho_{\rm eq}(x)+\begin{cases}10^{-2},&x\in[1,2],\\0,&\mbox{otherwise},\end{cases}\qquad q(x,0)=\sigma(x)\rho(x,0)u_{\rm eq}(x),
\end{equation*}
prescribed in the computational domain $[0,20]$ subject to free boundary conditions.

\begin{figure}[ht!]
\centerline{\includegraphics[trim=1.1cm 0.3cm 1.0cm 0.1cm, clip, width=4.cm]{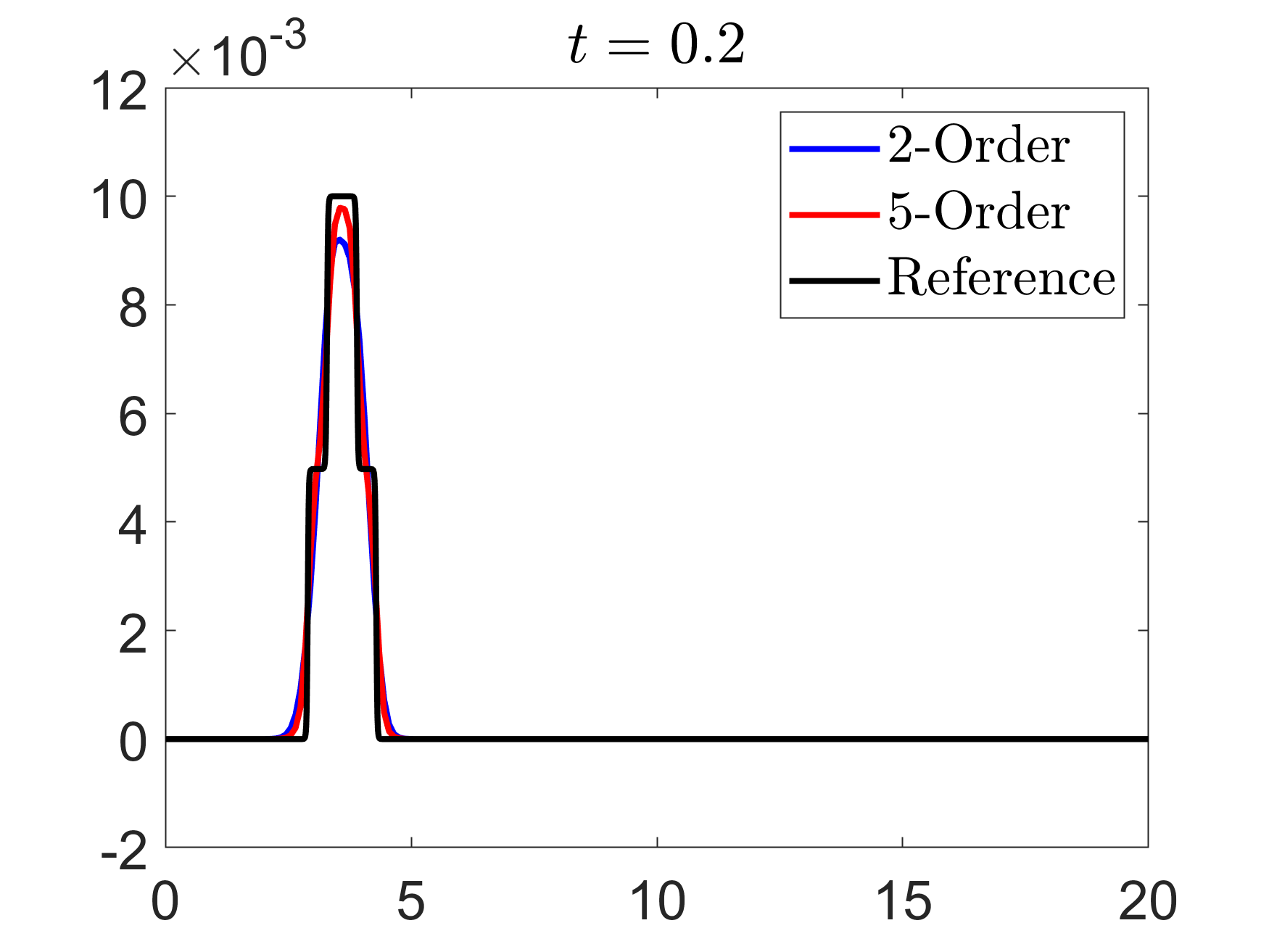}\hspace{0.1cm}
            \includegraphics[trim=1.1cm 0.3cm 1.0cm 0.1cm, clip, width=4.cm]{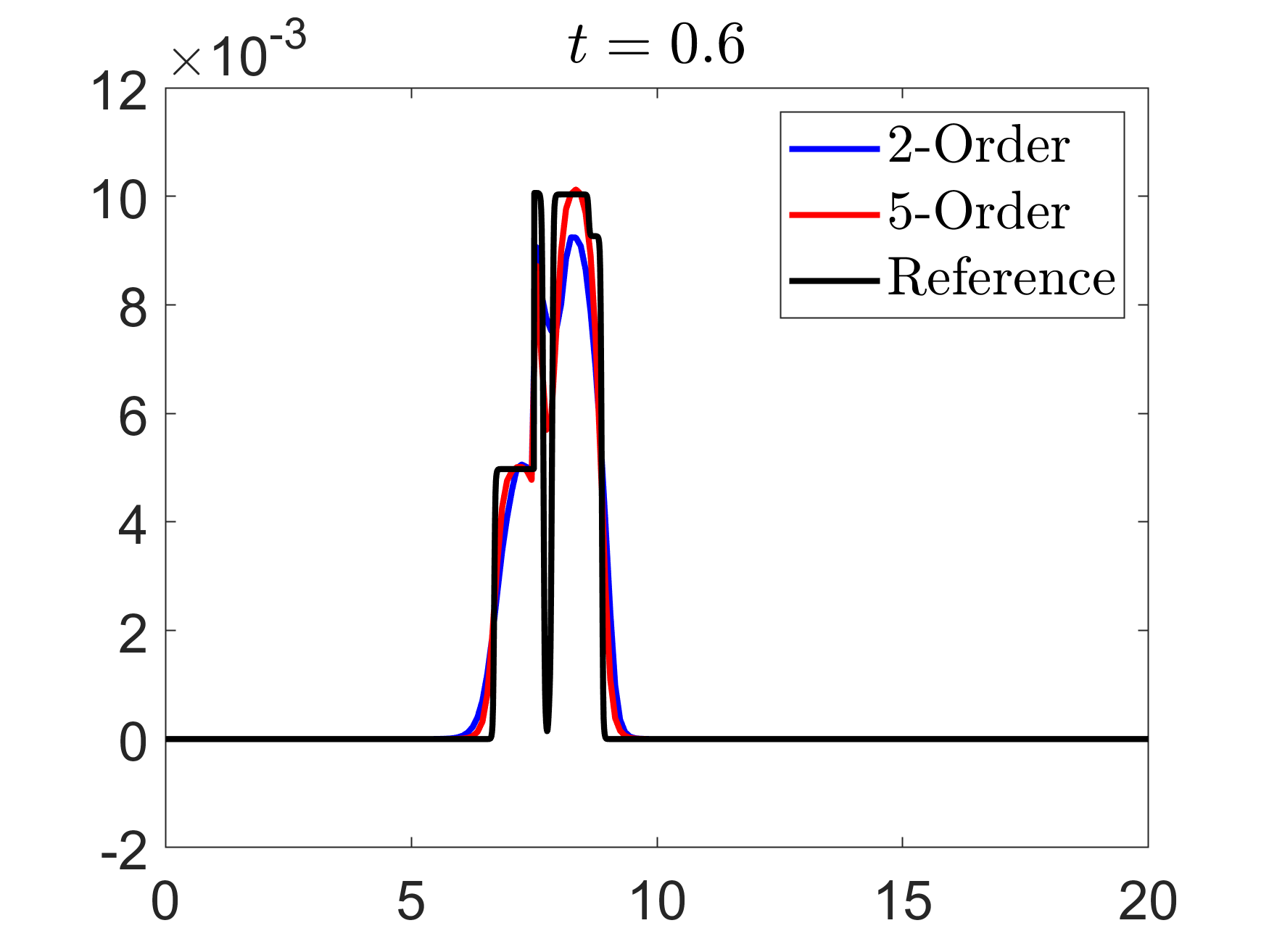}\hspace{0.1cm}
            \includegraphics[trim=1.1cm 0.3cm 1.0cm 0.1cm, clip, width=4.cm]{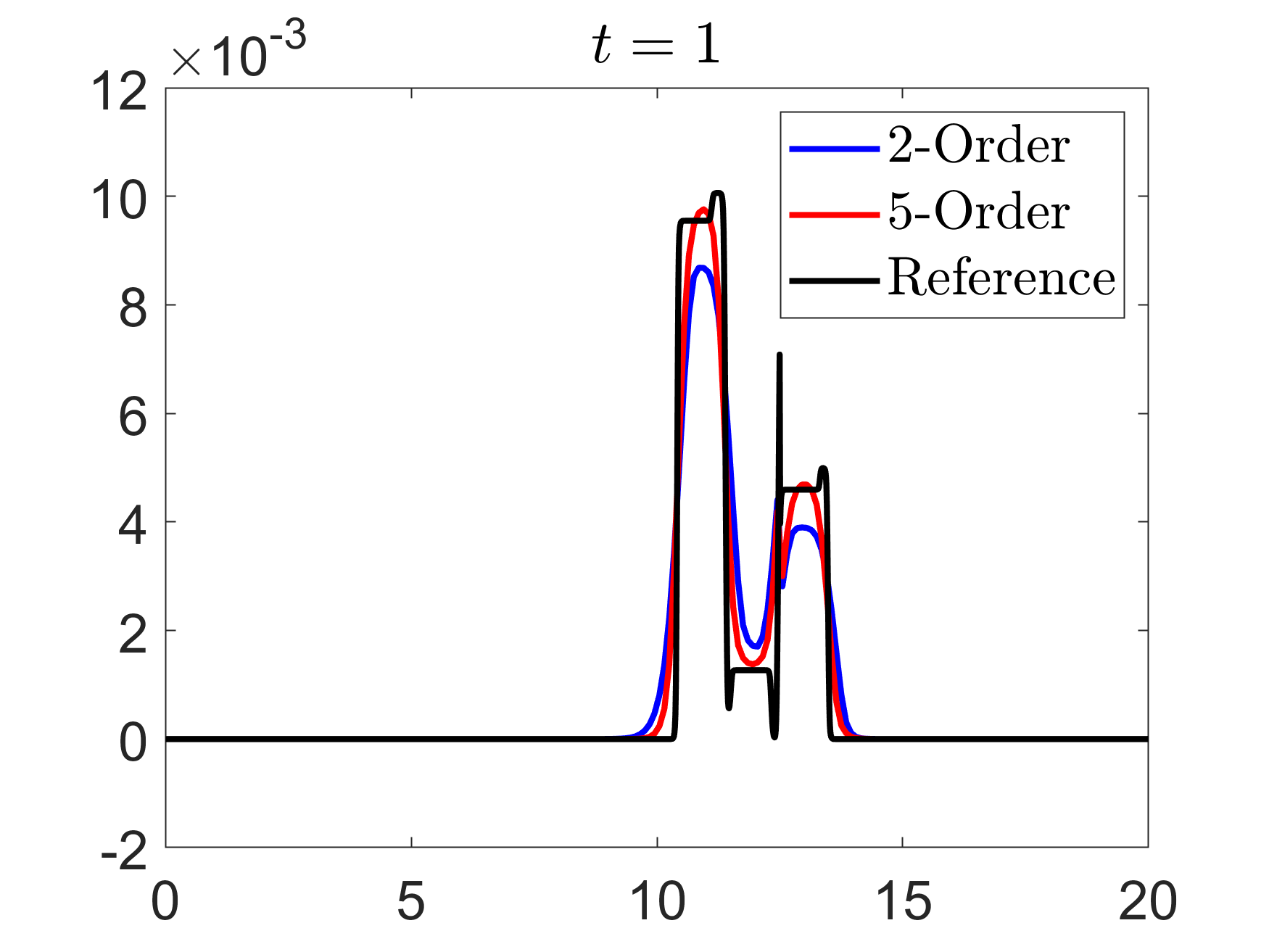}}
\vskip5pt
\centerline{\includegraphics[trim=1.1cm 0.3cm 1.0cm 0.1cm, clip, width=4.cm]{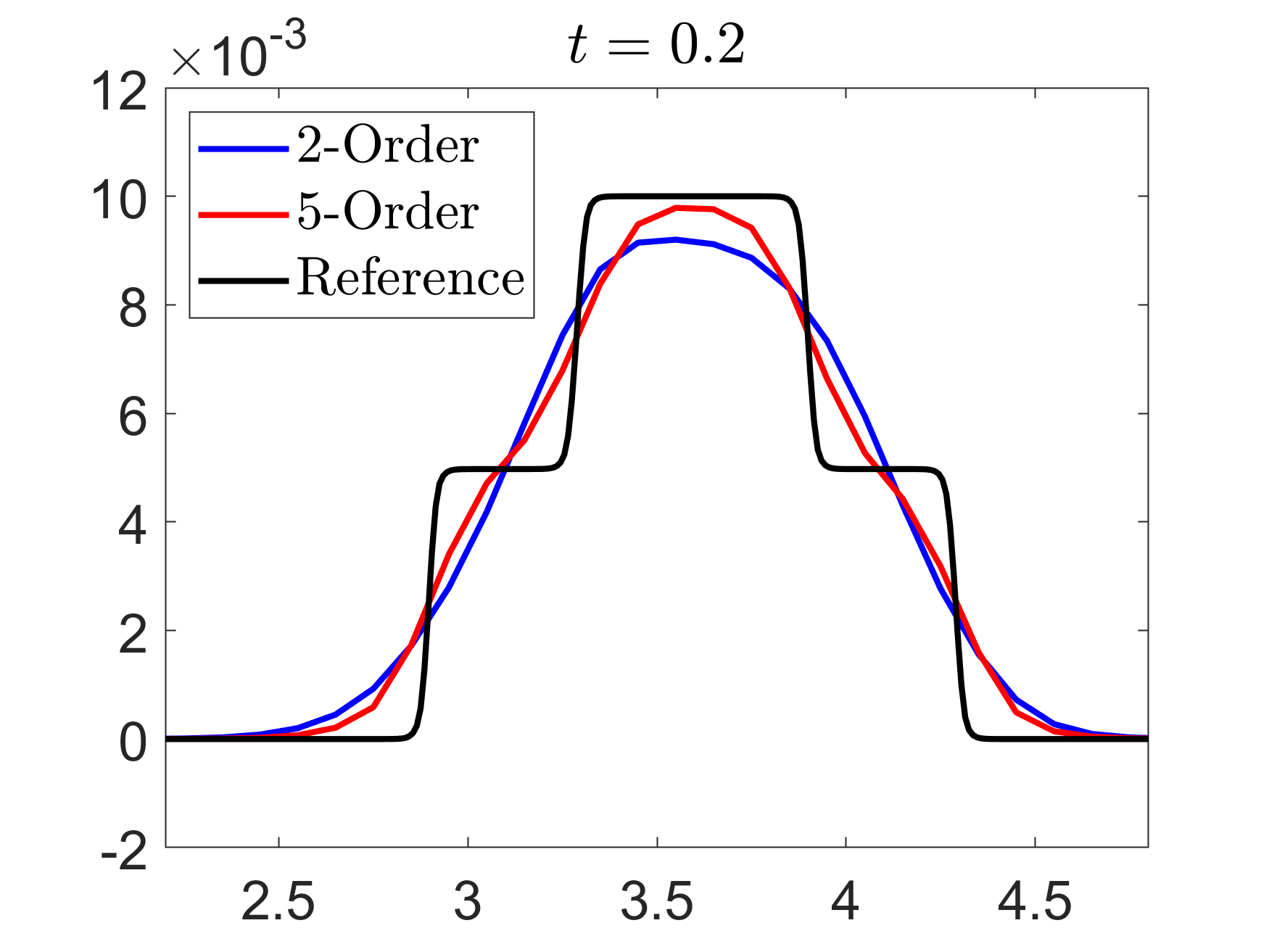}\hspace{0.1cm}
            \includegraphics[trim=1.1cm 0.3cm 1.0cm 0.1cm, clip, width=4.cm]{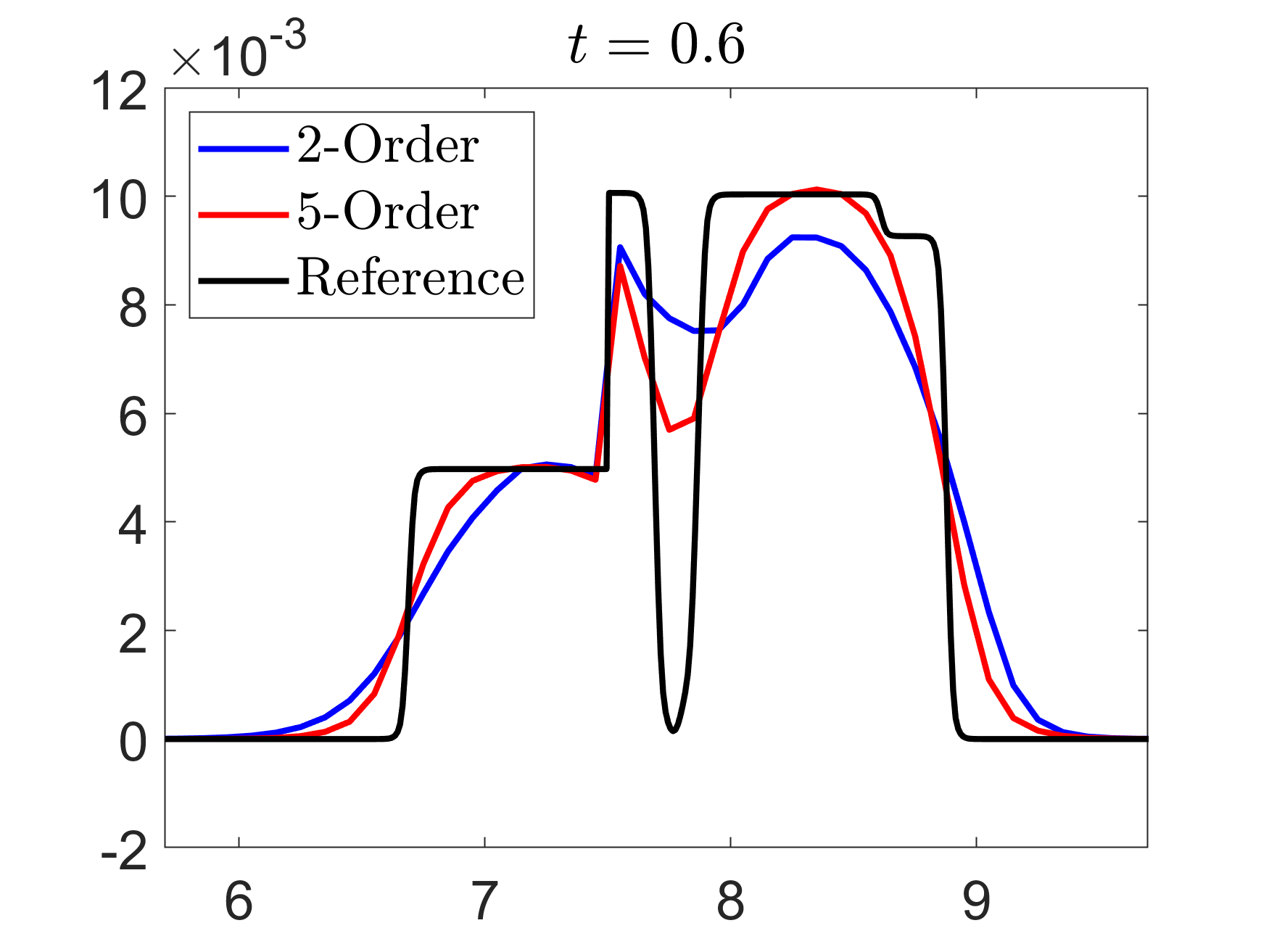}\hspace{0.1cm}
            \includegraphics[trim=1.1cm 0.3cm 1.0cm 0.1cm, clip, width=4.cm]{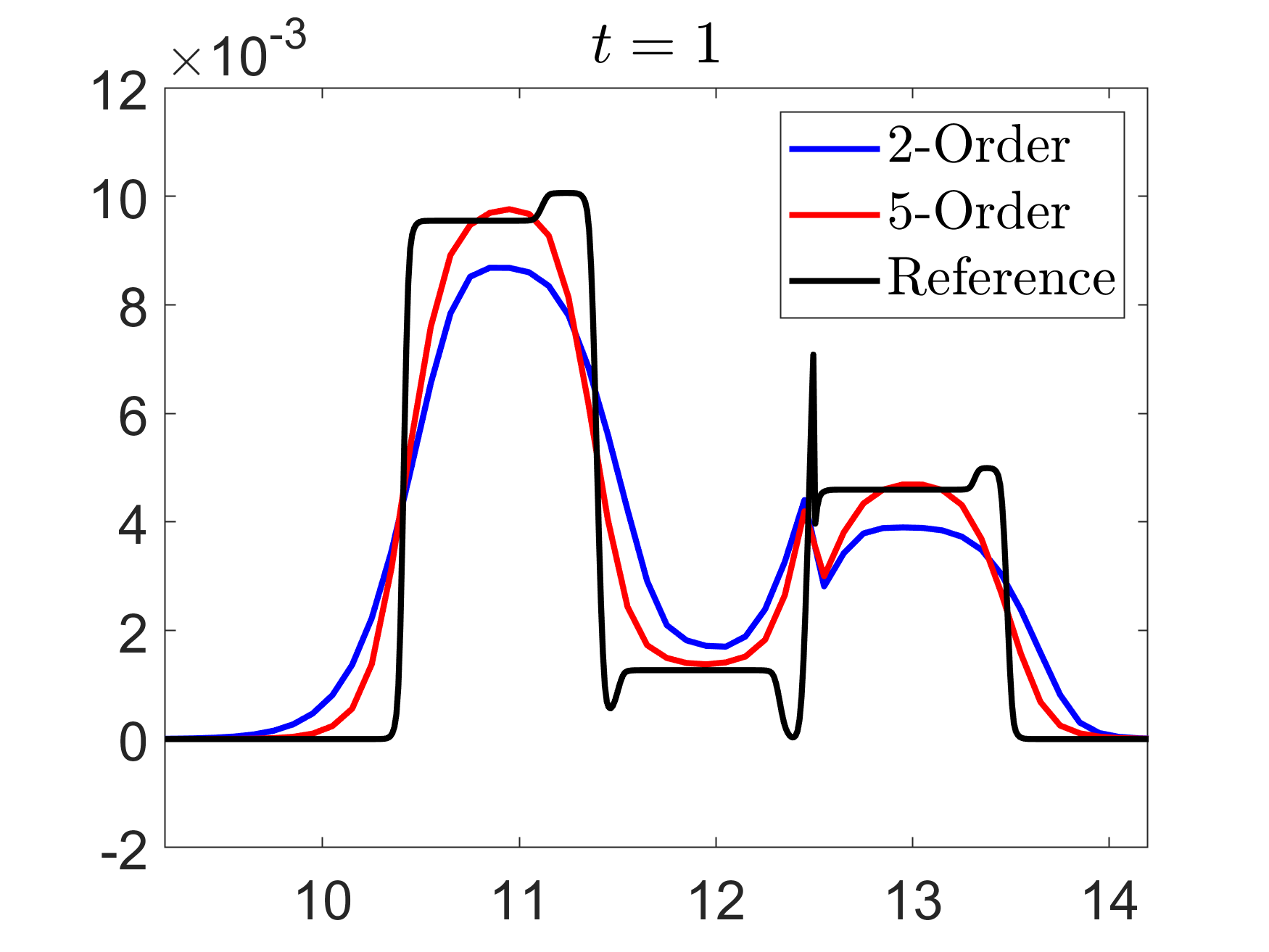}}
\caption{\sf Example 2: The difference $\rho(x,t)-\rho_{eq}$ (top row) and zoom at the perturbations (bottom row) at times $t=0.2$ (left
column), 0.6 (middle column), and 1 (right column).\label{fig9}}
\end{figure}

We compute the solutions until the final time $t=1$ by both the 2-Order and 5-Order Schemes on a uniform mesh with $\dx=1/10$. In Figure
\ref{fig9}, we plot the difference $\rho(x,t)-\rho_{\rm eq}(x)$ at times $t=0.2$, 0.6, and 1 together the reference solutions computed by
the 5-Order Scheme with $\dx=1/100$. The obtained results clearly demonstrate that the 5-Order Scheme achieves higher
resolution than the 2-Order Scheme, especially at larger times.

\subsection{Two-Layer Shallow Water Equations}
In this section, we proceed with three numerical examples for the two-layer shallow water equations \eref{1.1}, \eref{3.2.1}. We take
$r=0.98$ and either $g=10$ (Examples 3 and 5) or $g=1$ (Example 4).

\subsubsection*{Example 3---Experimental Order of Accuracy}
In this example taken from \cite{KPmultil}, we consider the following initial data:
\begin{equation*}
\begin{aligned}
&h_1(x,0)=5+e^{\cos(2\pi x)},\quad h_2(x,0)=5-e^{\cos(2\pi x)}-\sin^2(\pi x),\\
&q_1(x,0)=q_2(x,0)\equiv0,
\end{aligned}
\end{equation*}
and a continuous bottom topography
\begin{equation*}
Z(x)=\sin^2(\pi x)-10,
\end{equation*}
both prescribed in the computational domain $[0,1]$ subject to the periodic boundary conditions. 

We compute the numerical solution until the final time $t=0.1$ by both the 2-Order and 5-Order Schemes on a sequence of uniform meshes with
$\dx=1/40$, 1/80, 1/160, 1/320, 1/640, and 1/1280. We measure the $L^1$-errors in $h_1$ and estimate the experimental convergence rates
using the following Runge formulae, which are based on the solutions computed on the three consecutive uniform grids with the mesh sizes
$\dx$, $2\dx$, and $4\dx$ and denoted by $(\cdot)^{\dx}$, $(\cdot)^{2\dx}$, and $(\cdot)^{4\dx}$, respectively:
$$
{\rm Error}(\dx)\approx\frac{\delta_{12}^2}{|\delta_{12}-\delta_{24}|},\quad
{\rm Rate}(\dx)\approx\log_2\left(\frac{\delta_{24}}{\delta_{12}}\right).
$$
Here, $\delta_{12}:=\|(\cdot)^{\dx}-(\cdot)^{2\dx}\|_{L^1}$ and $\delta_{24}:=\|(\cdot)^{2\dx}-(\cdot)^{4\dx}\|_{L^1}$. The obtained results
for the density are reported in Table \ref{tab71}, where one can clearly see that the expected orders of accuracy are achieved by both the
2-Order and 5-Order Schemes. Note that in order to achieve the fifth order of accuracy, we have used smaller time steps with
$\dt\sim(\dx)^\frac{5}{3}$.
\begin{table}[ht!]
\caption{\sf Example 3: The $L^1$-errors and experimental convergence rates for $h_1$.\label{tab71}}
\centering
\begin{tabular}{|c|cc|cc|cc|cc|}
\hline
\multirow{2}{1em}{$\dx$}&\multicolumn{2}{c|}{2-Order}&\multicolumn{2}{c|}{5-Order}\\
\cline{2-5}&Error&Rate&Error&Rate\\
\hline
$1/160$ & 5.64e-4 &2.58&1.90e-7 &4.83\\
$1/320$ & 6.18e-5 &2.86&5.10e-9 &5.02\\
$1/640$ & 1.38e-5 &2.54&1.57e-10&5.02\\
$1/1280$& 1.75e-6 &2.74&4.86e-12&5.01\\
\hline
\end{tabular}
\end{table}

\subsubsection*{Example 4---Small Perturbation of Discontinuous Steady State}
In this example, we consider a discontinuous steady state given by 
\begin{equation}
\begin{aligned}
(h_1)_{\rm eq}(x):&=\begin{cases}1.22373355048230,&x<0,\\1.44970064153589,&x>0,\end{cases}\qquad(q_1)_{\rm eq}(x)\equiv12,\\
(h_2)_{\rm eq}(x):&=\begin{cases}0.968329515483846,&x<0,\\1.12439026921484,&x>0,\end{cases}\qquad(q_2)_{\rm eq}(x)\equiv10,\\
\end{aligned}
\label{20}
\end{equation}
and a discontinuous bottom topography
\begin{equation}
Z(x)=\begin{cases}-2,&x<0,\\-1,&x>0.\end{cases}
\label{21}
\end{equation}

In order to demonstrate that the proposed 2-Order and 5-Order Schemes are WB, we use \eref{20}--\eref{21} as the initial setting
prescribed in the computational domain $[-1,1]$ subject to free boundary conditions and compute the numerical solution until the final time
$t=20$ on a uniform mesh with $\dx=1/100$. The obtained discrete $L^1$- and $L^\infty$-errors are reported in Table \ref{tab1}, where one
can clearly see that both the 2-Order and 5-Order Schemes can preserve the steady state within the machine accuracy.
\begin{table}[ht!]
\caption{\sf Example 4: $L^1$-errors for $h_1$, $h_2$, $q_1$, and $q_2$.\label{tab1}}
\centering
\begin{tabular}{ccccccc|cc|}
\hline
&$L^1$-error in $h_1$ &$L^1$-error in $h_2$&$L^1$-error in $q_1$&$L^1$-error in $q_2$\\
\hline
2-Order& 1.77e-17& 2.58e-17& 1.42e-16 & 2.13e-16\\
5-Order& 4.44e-17& 3.37e-16& 1.40e-15 &4.96e-15 \\
\hline
\end{tabular}
\end{table}

We then test the ability of the studied 2-Order and 5-Order Schemes to capture small perturbations of this steady state, which is initially
added to the upper layer depth:
\begin{equation*}
\begin{aligned}
h_1(x,0)&=(h_1)_{\rm eq}(x)+\begin{cases}0.001,&x\in[-0.6,-0.5],\\0,&\mbox{otherwise},\end{cases}\\
h_2(x,0)&=(h_2)_{\rm eq}(x), \quad q_1(x,0)=(q_1)_{\rm eq}(x),\quad q_2(x,0)=(q_2)_{\rm eq}(x).
\end{aligned}
\end{equation*}

\begin{figure}[ht!]
\centerline{\includegraphics[trim=0.0cm 0.3cm 0.7cm 0.2cm, clip, width=4.cm]{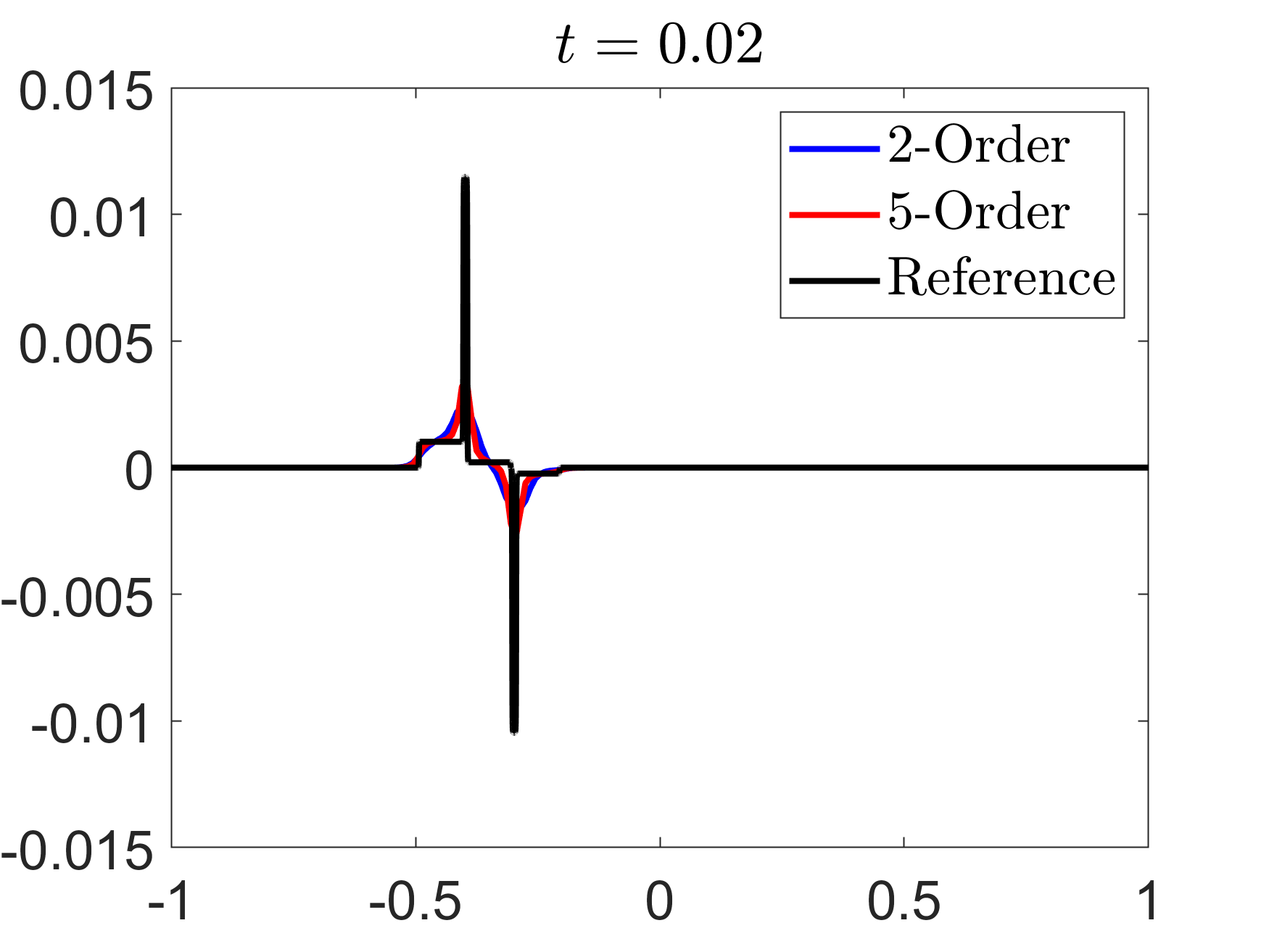}\hspace{0.1cm}
            \includegraphics[trim=0.0cm 0.3cm 0.7cm 0.2cm, clip, width=4.cm]{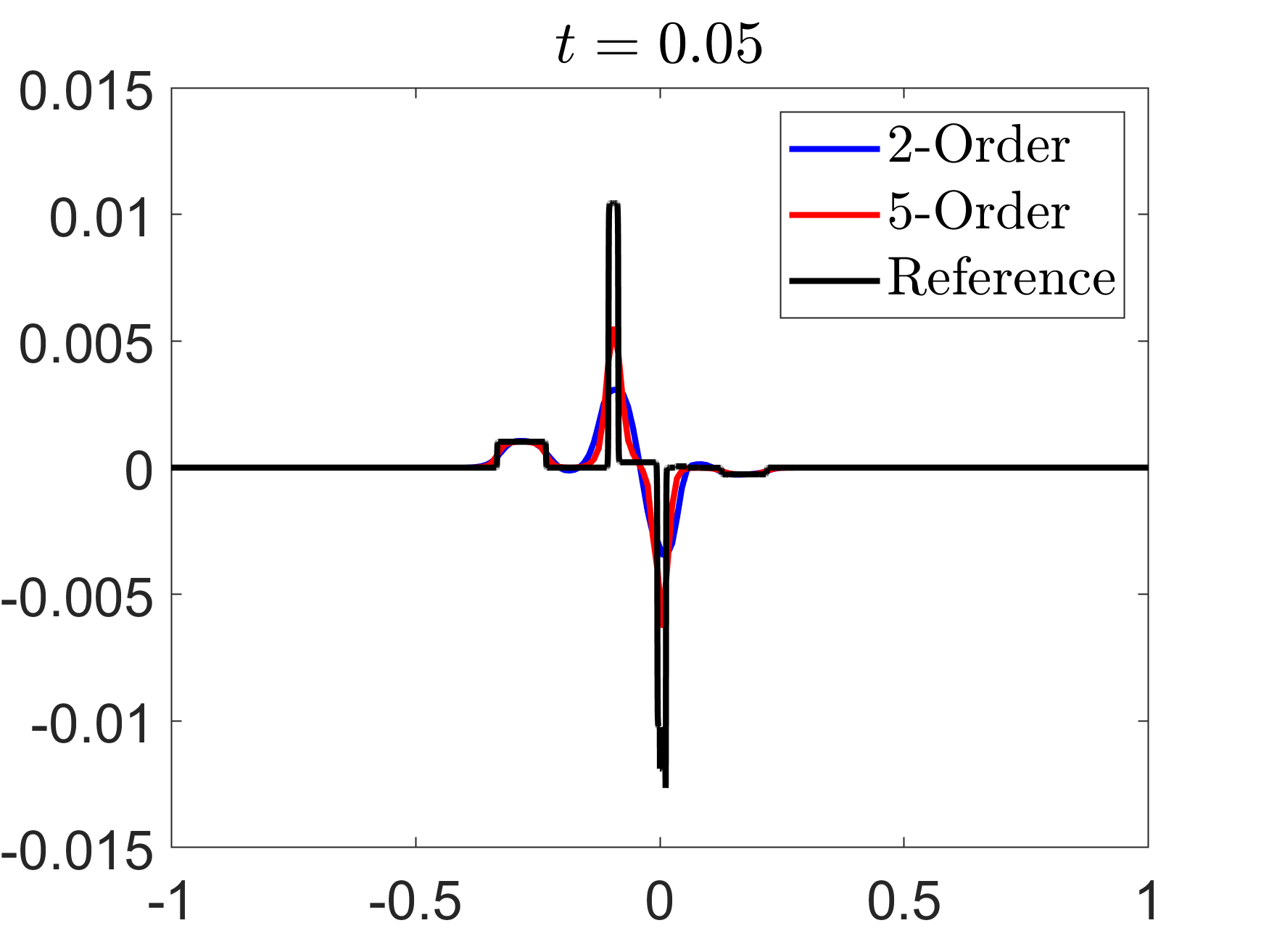}\hspace{0.1cm}
            \includegraphics[trim=0.0cm 0.3cm 0.7cm 0.2cm, clip, width=4.cm]{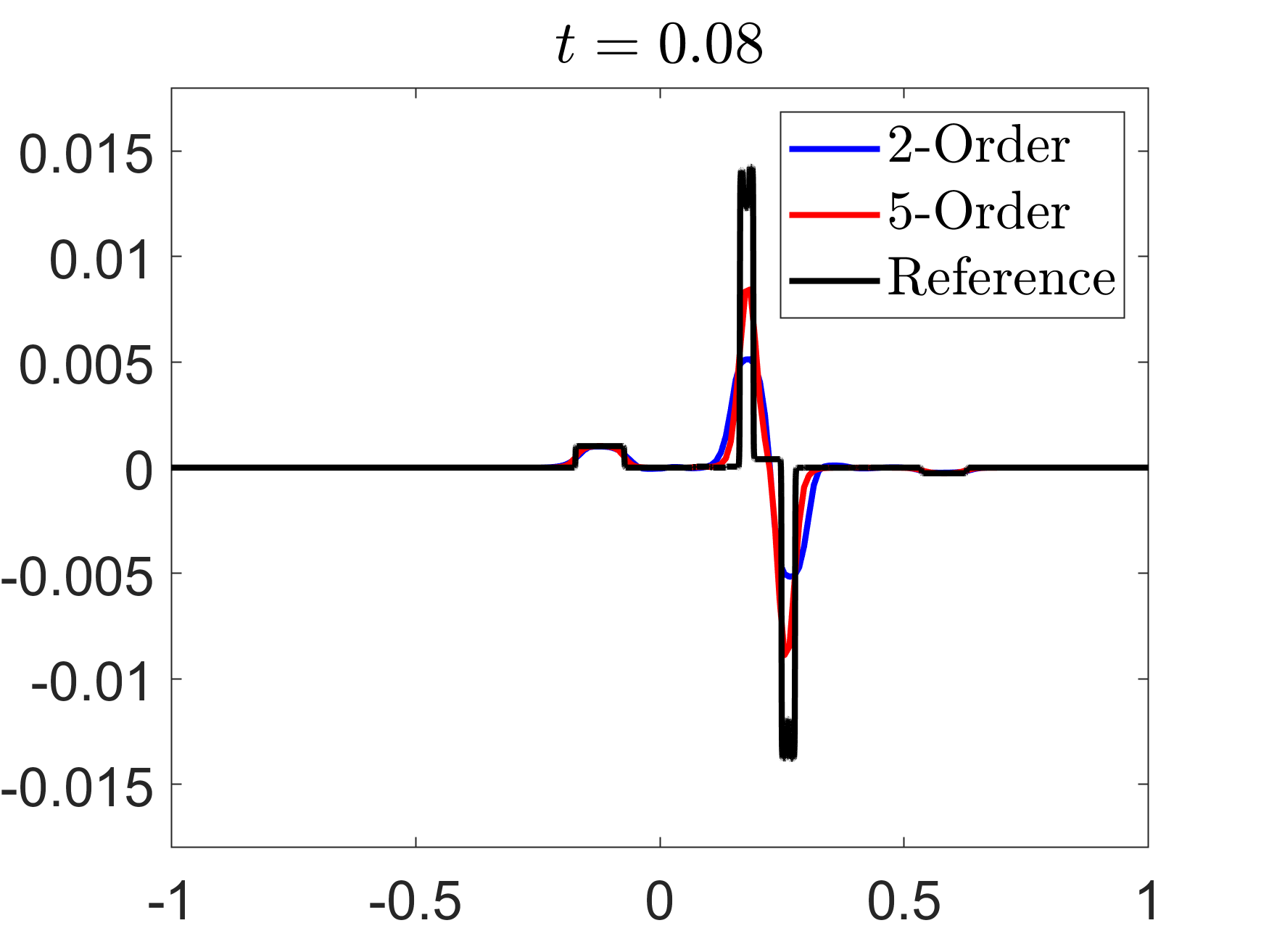}}
\vskip5pt
\centerline{\includegraphics[trim=0.0cm 0.3cm 0.7cm 0.2cm, clip, width=4.cm]{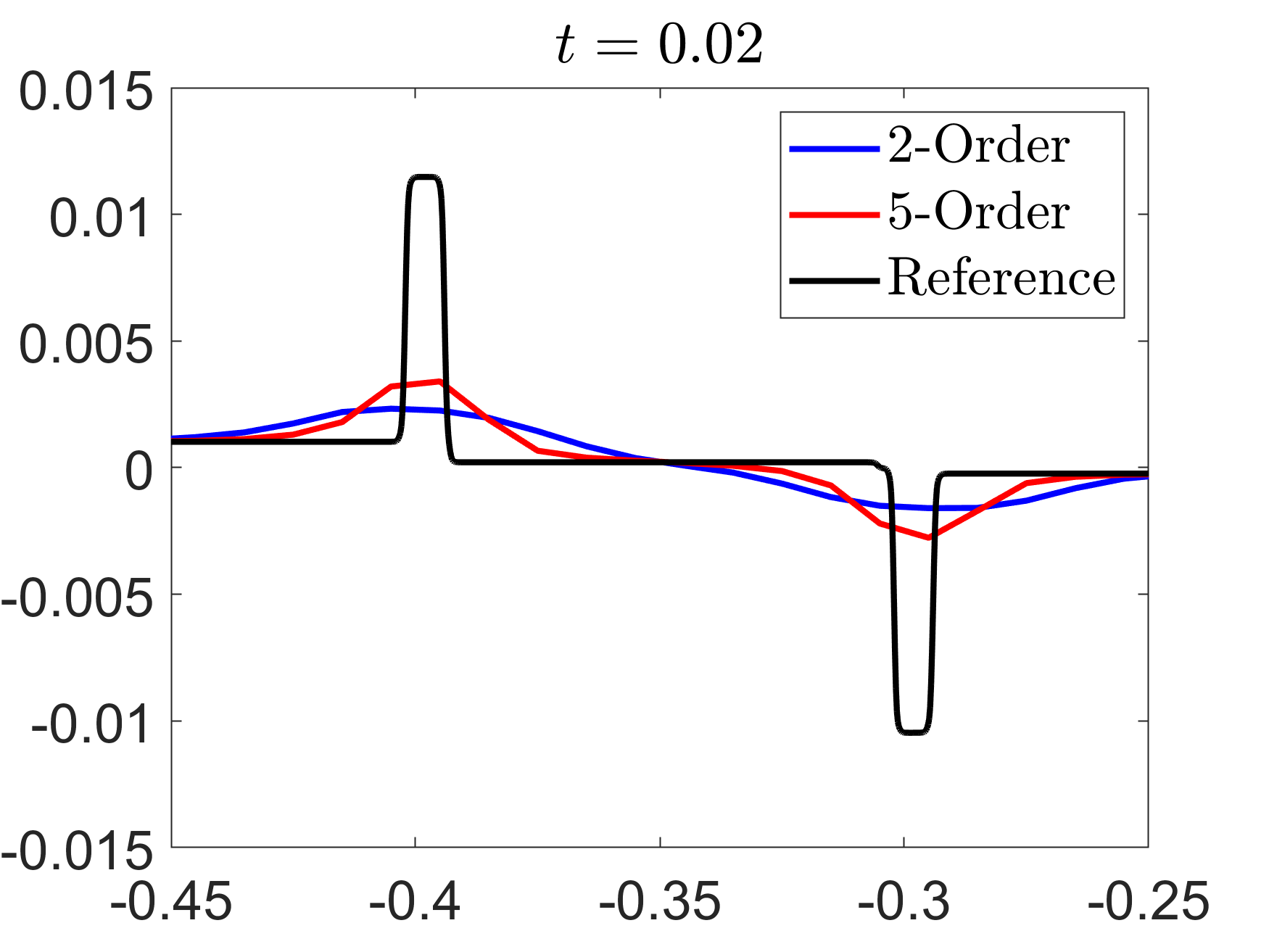}\hspace{0.1cm}
            \includegraphics[trim=0.0cm 0.3cm 0.7cm 0.2cm, clip, width=4.cm]{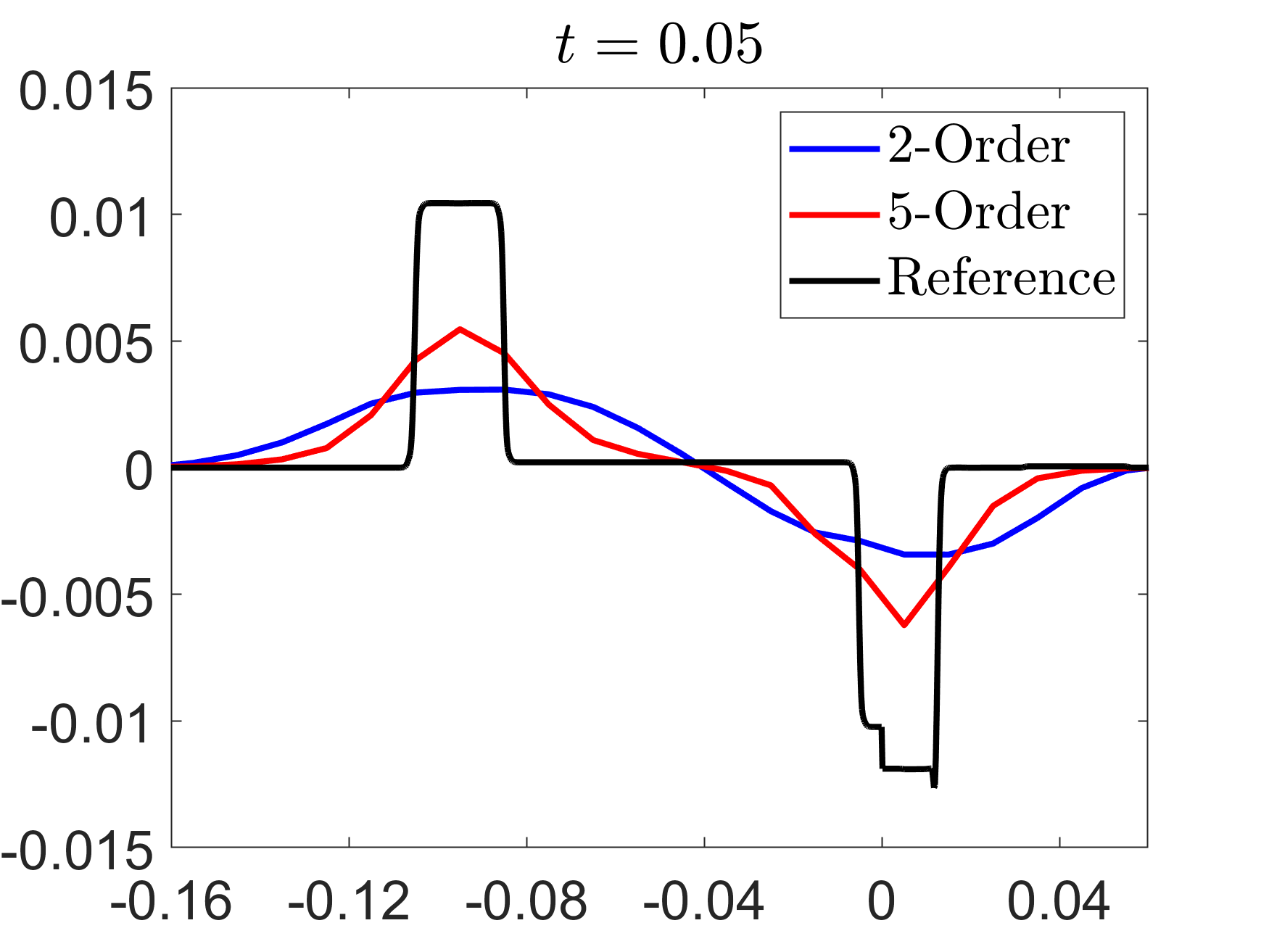}\hspace{0.1cm}
            \includegraphics[trim=0.0cm 0.3cm 0.7cm 0.2cm, clip, width=4.cm]{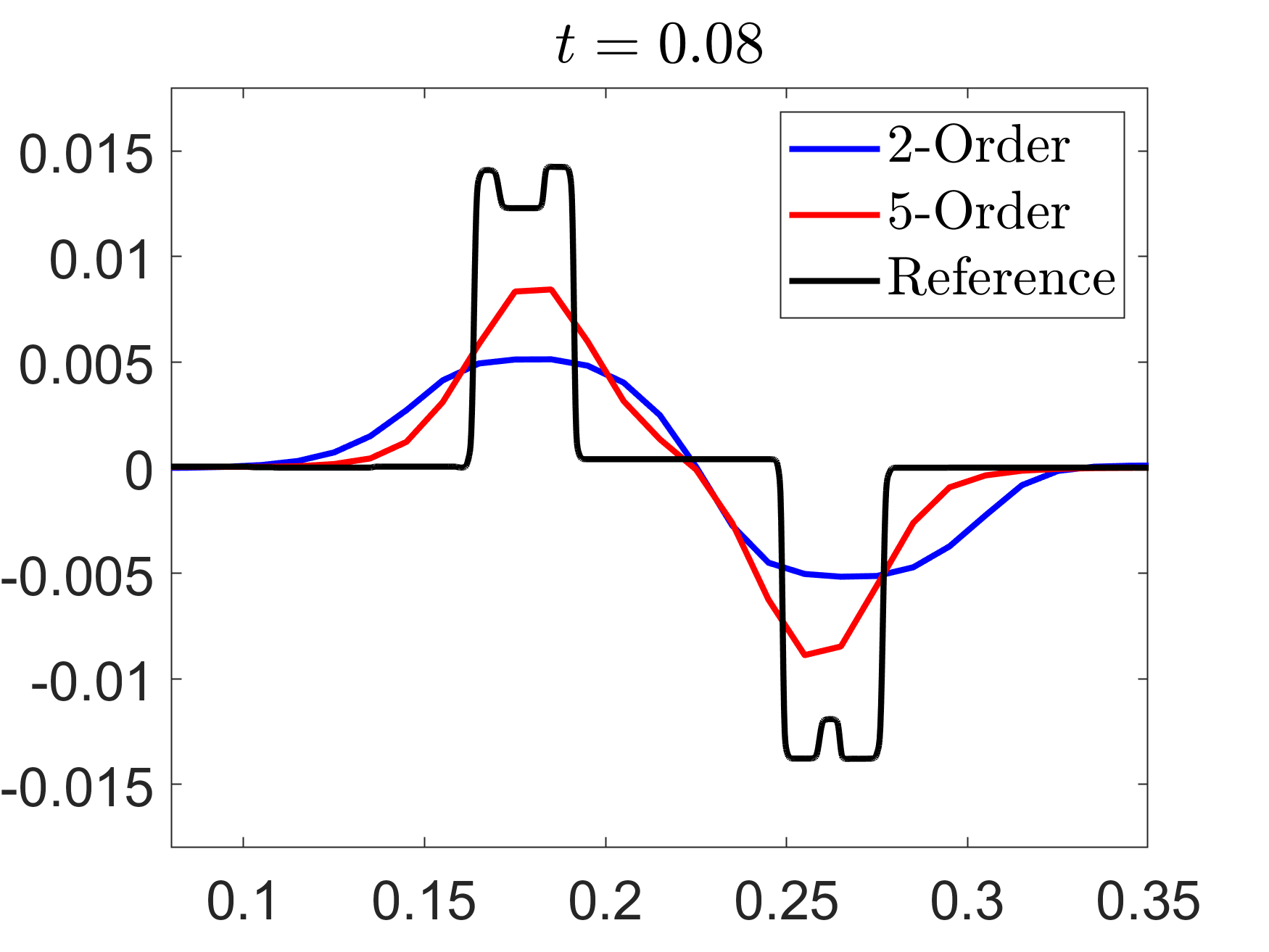}}
\caption{\sf Example 4: The difference $h_1(x,t)-(h_1)_{\rm eq}$ (top row) and zoom at the perturbations (bottom row) at times $t=0.02$
(left column), 0.05 (middle column), and 0.08 (right column).\label{fig11}}
\end{figure}

We compute the numerical solution until the final time $t=0.08$ by both the 2-Order and 5-Order Schemes on a uniform mesh with $\dx=1/100$
and plot the difference $h_1(x,t)-(h_1)_{\rm eq}(x)$ at times $t=0.02$, 0.05, and 0.08 in Figure \ref{fig11}. Here, the reference solution
is computed by the 5-Order Scheme on a much finer mesh with $\dx=1/5000$. As one can see, the 5-Order Scheme outperforms the
2-Order Scheme.

\subsubsection*{Example 5---Convergence to a Moving-Water Steady State}
In the fourth example, we first demonstrate the ability of the proposed 5-Order Scheme to reach discrete steady state. To this end, we take
constant initial data
\begin{equation*}
h_1(x,0)\equiv8,\quad h_2(x,0)\equiv4,\quad q_1(x,0)=q_2(x,0)\equiv0,
\end{equation*}
and a discontinuous bottom topography
\begin{equation*}
Z(x)=\begin{cases}0.2,&x\in[8,12],\\0,&\mbox{otherwise},\end{cases}
\end{equation*}
both prescribed in the computational domain $[0,25]$ subject to the Dirichlet boundary conditions imposed at $x=0$,
\begin{equation*}
h_1(0,t)=8,\quad(h_2)(0,t)=4,\quad(q_1)(0,t)=119,\quad(q_2)(0,t)=60,
\end{equation*}
and free boundary conditions imposed at $x=25$.

We first compute the numerical solutions on a uniform mesh with $\dx=1/5$ until a very large time $t=100$ to let the solution to converge to
the discrete steady state. The obtained results agree almost perfectly with those reported in \cite[Figure 7.20]{KLX_21} for the 2-Order
Scheme.

We then add a small perturbation to the upper layer water depth and consider the following initial data:
\begin{equation*}
\begin{aligned}
h_1(x,0)&=(h_1)_{\rm eq}(x)+\begin{cases}0.0001,&x\in[2,2.25],\\0,&\mbox{otherwise},\end{cases}\\
h_2(x,0)&=(h_2)_{\rm eq}(x),\quad q_1(x,0)=(q_1)_{\rm eq}(x),\quad q_2(x,0)=(q_2)_{\rm eq}(x),
\end{aligned}
\end{equation*}
subject to free boundary conditions at both ends of the computational domain.

\begin{figure}[ht!]
\centerline{\includegraphics[trim=1.1cm 0.3cm 1.0cm 0.2cm, clip, width=4.0cm]{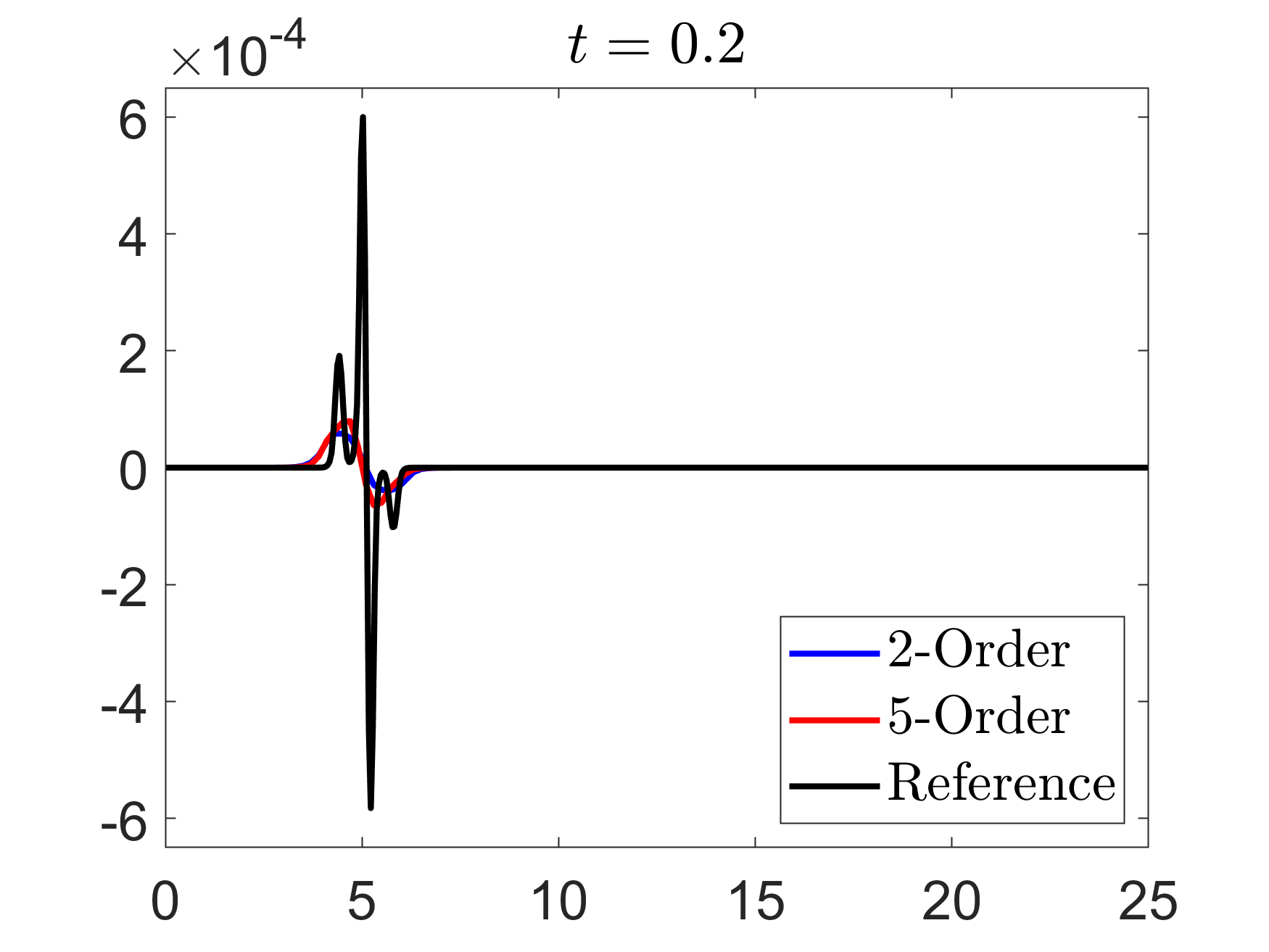}\hspace{0.1cm}
            \includegraphics[trim=1.1cm 0.3cm 1.0cm 0.2cm, clip, width=4.0cm]{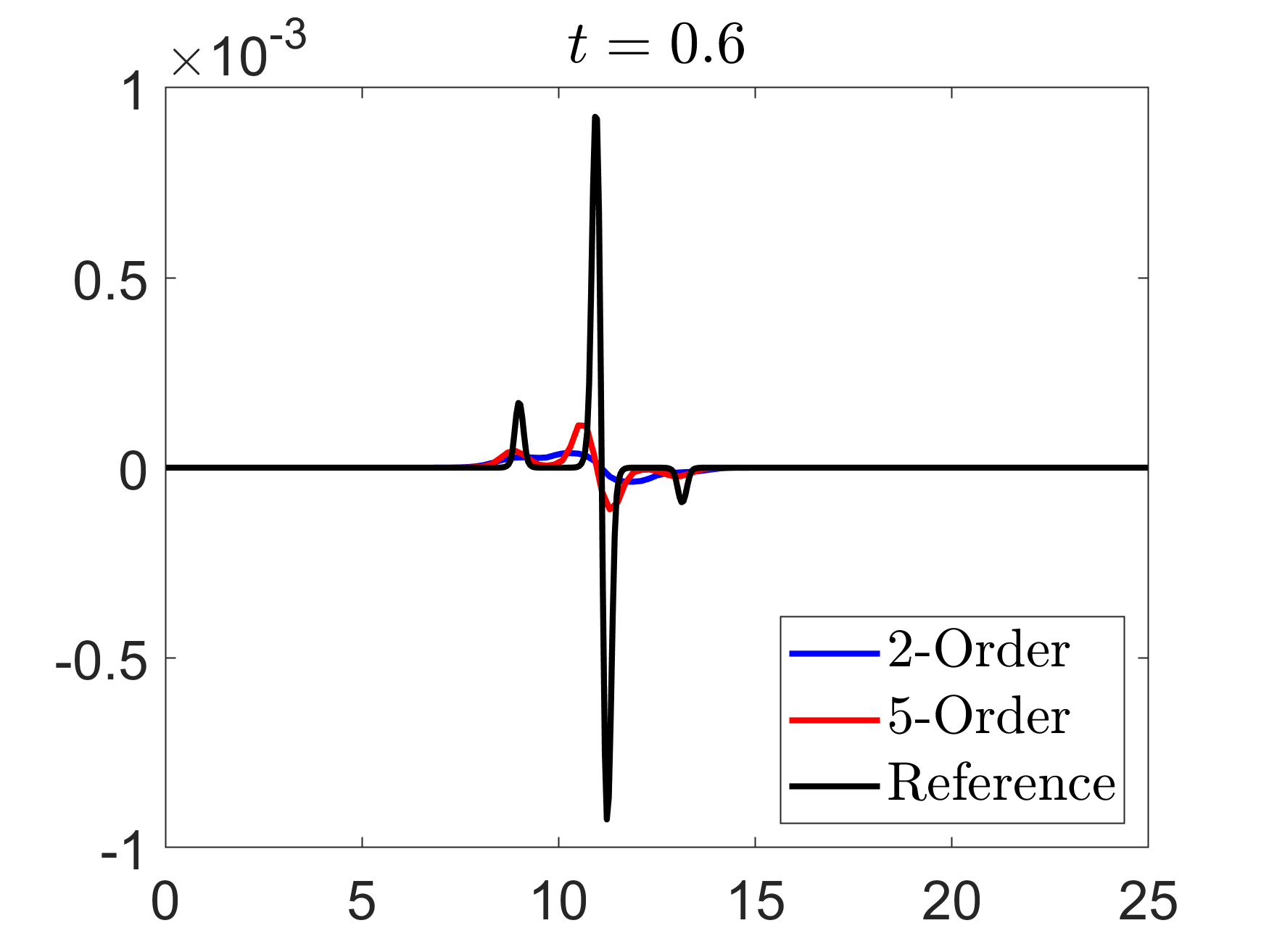}\hspace{0.1cm}
            \includegraphics[trim=1.1cm 0.3cm 1.0cm 0.2cm, clip, width=4.0cm]{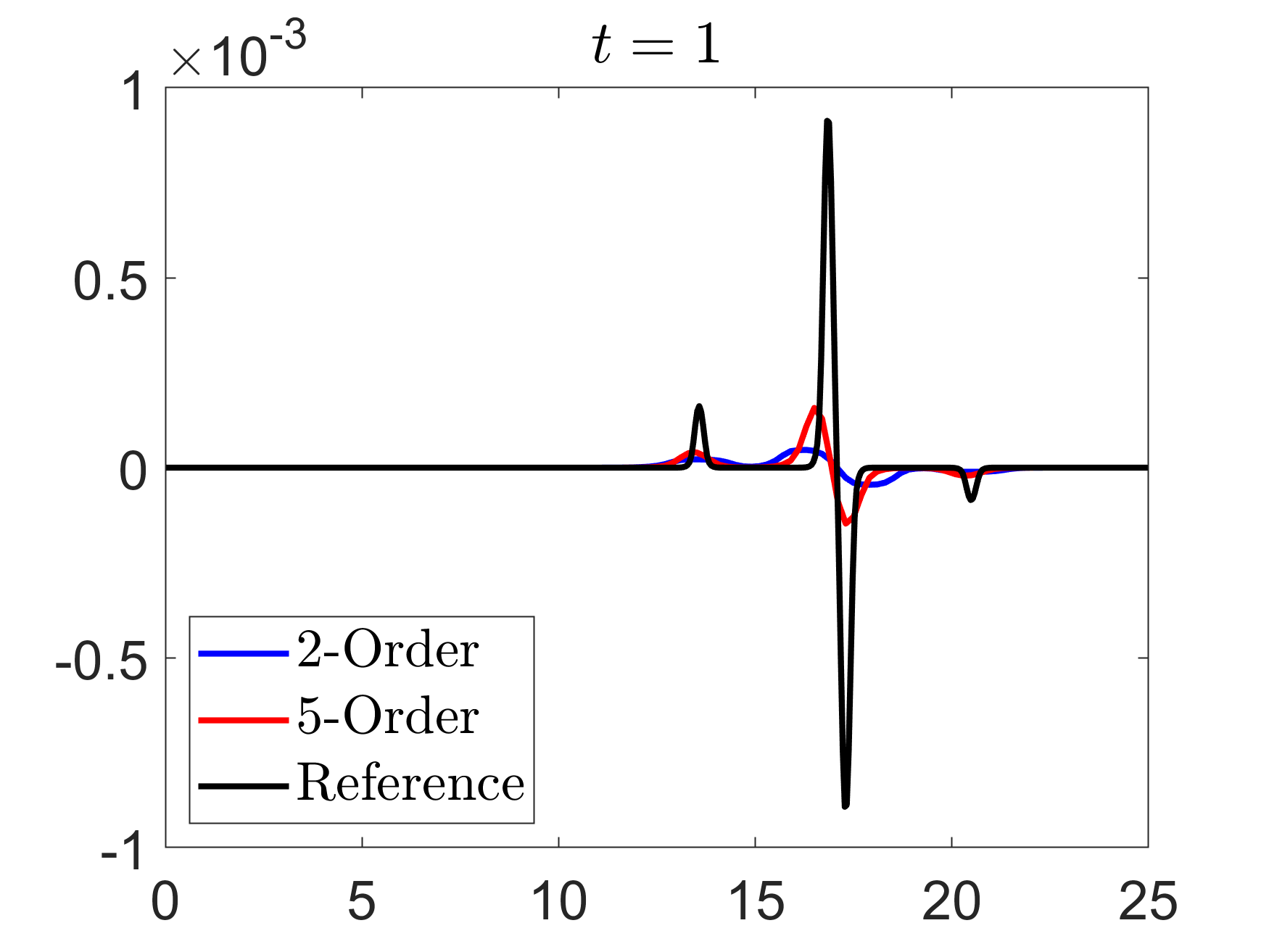}}
\vskip5pt
\centerline{\includegraphics[trim=1.1cm 0.3cm 1.0cm 0.2cm, clip, width=4.0cm]{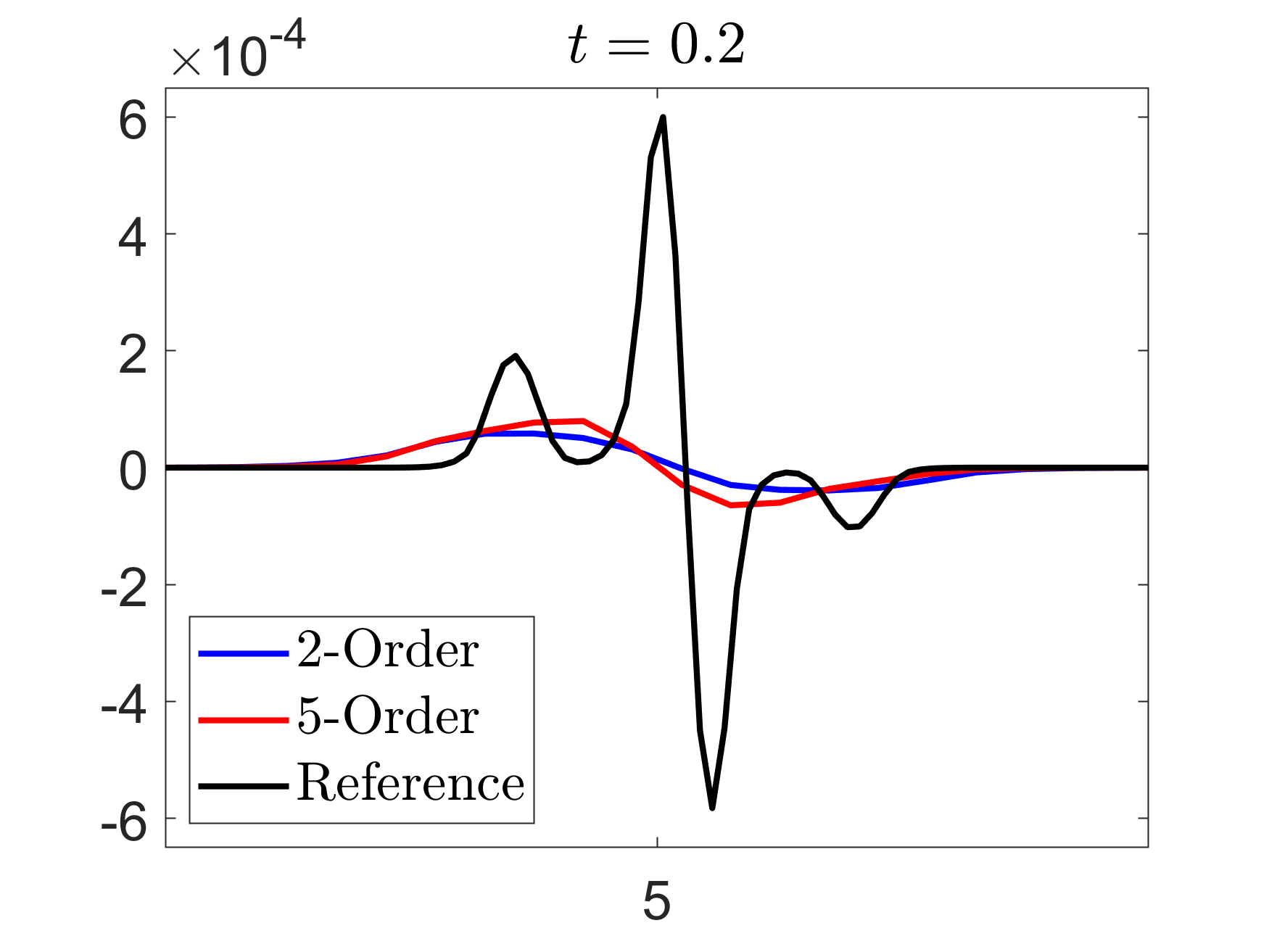}\hspace{0.1cm}
            \includegraphics[trim=1.1cm 0.3cm 1.0cm 0.2cm, clip, width=4.0cm]{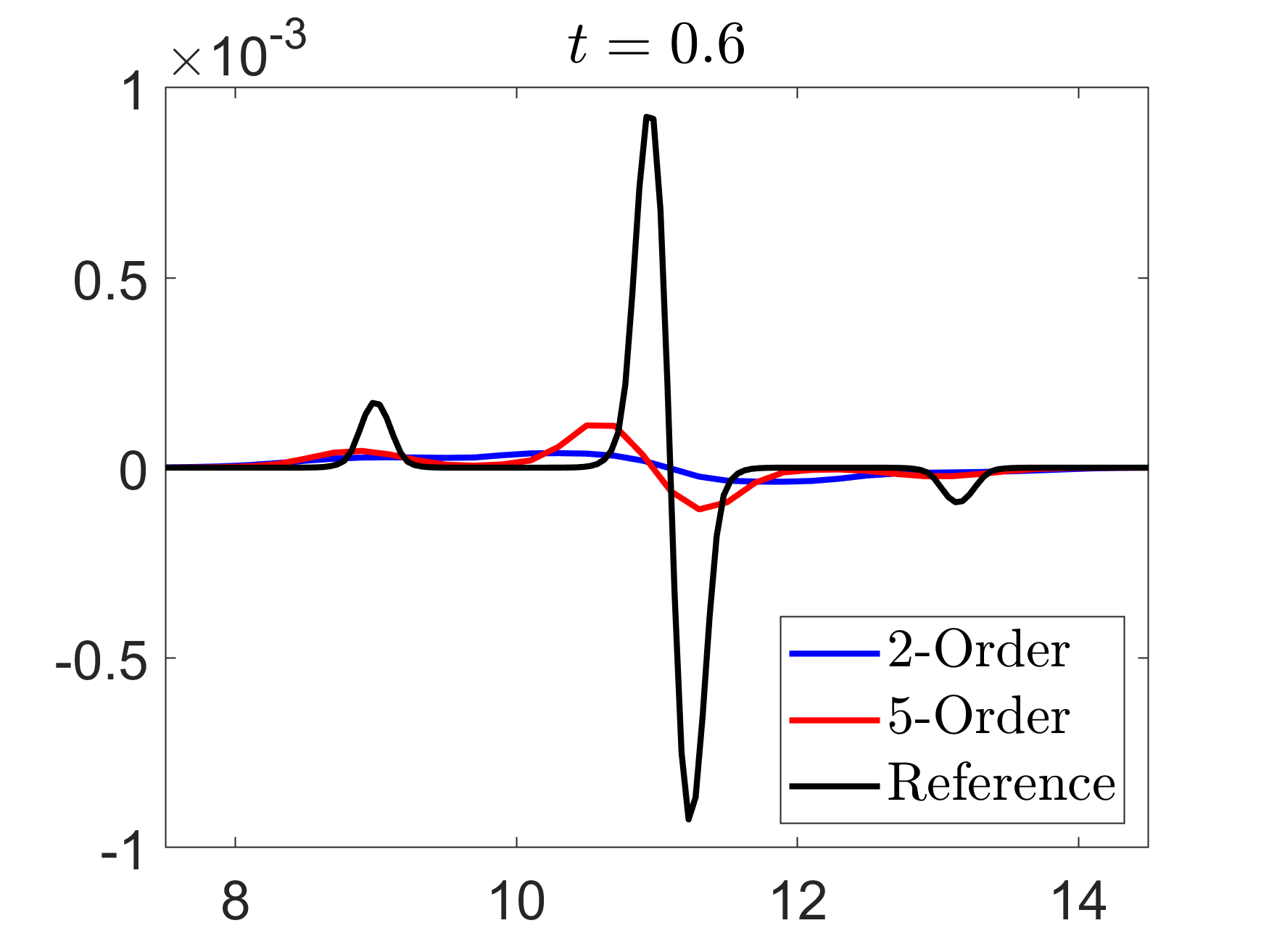}\hspace{0.1cm}
            \includegraphics[trim=1.1cm 0.3cm 1.0cm 0.2cm, clip, width=4.0cm]{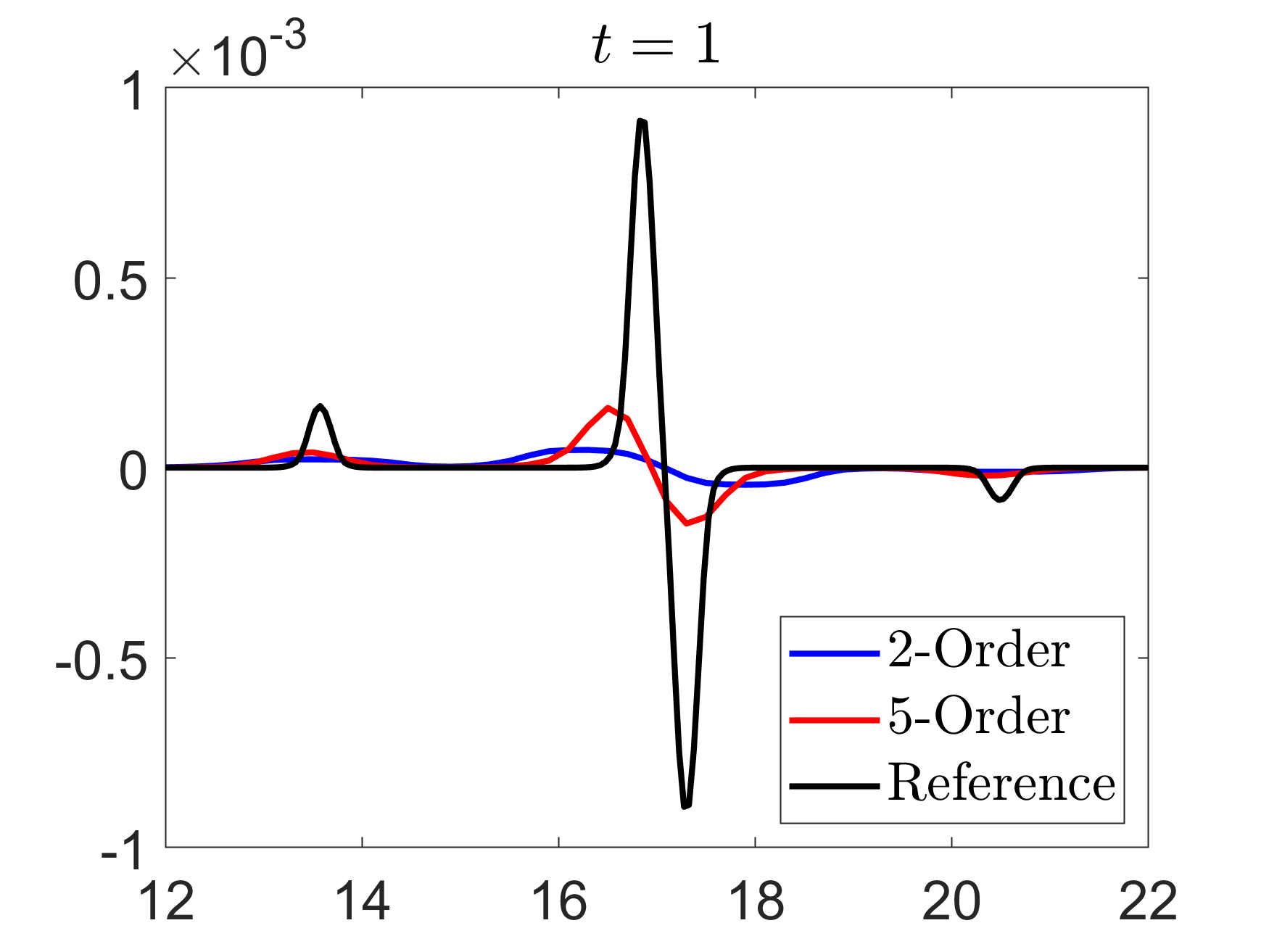}}
\caption{\sf Example 5: The difference $h_1(x,t)-(h_1)_{\rm eq}$ (top row) and zoom at the perturbations (bottom row) at times $t=0.2$ (left
column), 0.6 (middle column), and 1 (right column).\label{fig14}}
\end{figure}

We compute the numerical solutions until the final time $t=1$ by both the 2-Order and 5-Order Schemes on a uniform mesh with $\dx=1/5$. The
obtained differences $h_1-(h_1)_{\rm eq}$ at times $t=0.2$, 0.6, and 1 are plotted in Figure \ref{fig14} together with the reference
solution computed by the 5-Order Scheme on a finer mesh with $\dx=1/20$. One can observe that the 5-Order Scheme produces sharper results
compared to those obtained by the 2-Order Scheme.

\subsubsection*{Example 6---Riemann Problems}
In the final example, we numerically solve two Riemann problems with the following initial data:
\begin{equation*}
\begin{aligned}
{\rm\bf Test 1}:\quad(h_1,q_1,h_2,q_2)(x,0)&=\begin{cases}(1,1.5,1,1),&x<0,\\(0.8,1.2,1.2,1.8),&\mbox{otherwise},\end{cases}\\
{\rm\bf Test 2}:\quad(h_1,q_1,h_2,q_2)(x,0)&=\begin{cases}(1.5,1,1,1.5),&x<0,\\(1.2,1.6,0.9,1.2),&\mbox{otherwise},\end{cases}
\end{aligned}
\end{equation*}
with a discontinuous bottom topography,
\begin{equation*}
Z(x)=\begin{cases}-2,&x<0,\\-1.5,&\mbox{otherwise},\end{cases}
\end{equation*}
prescribed in the computational domain $[-1,1]$ subject to free boundary conditions.

\begin{figure}[ht!]
\centerline{\includegraphics[trim=0.8cm 0.3cm 1.cm 0.2cm, clip, width=5.cm]{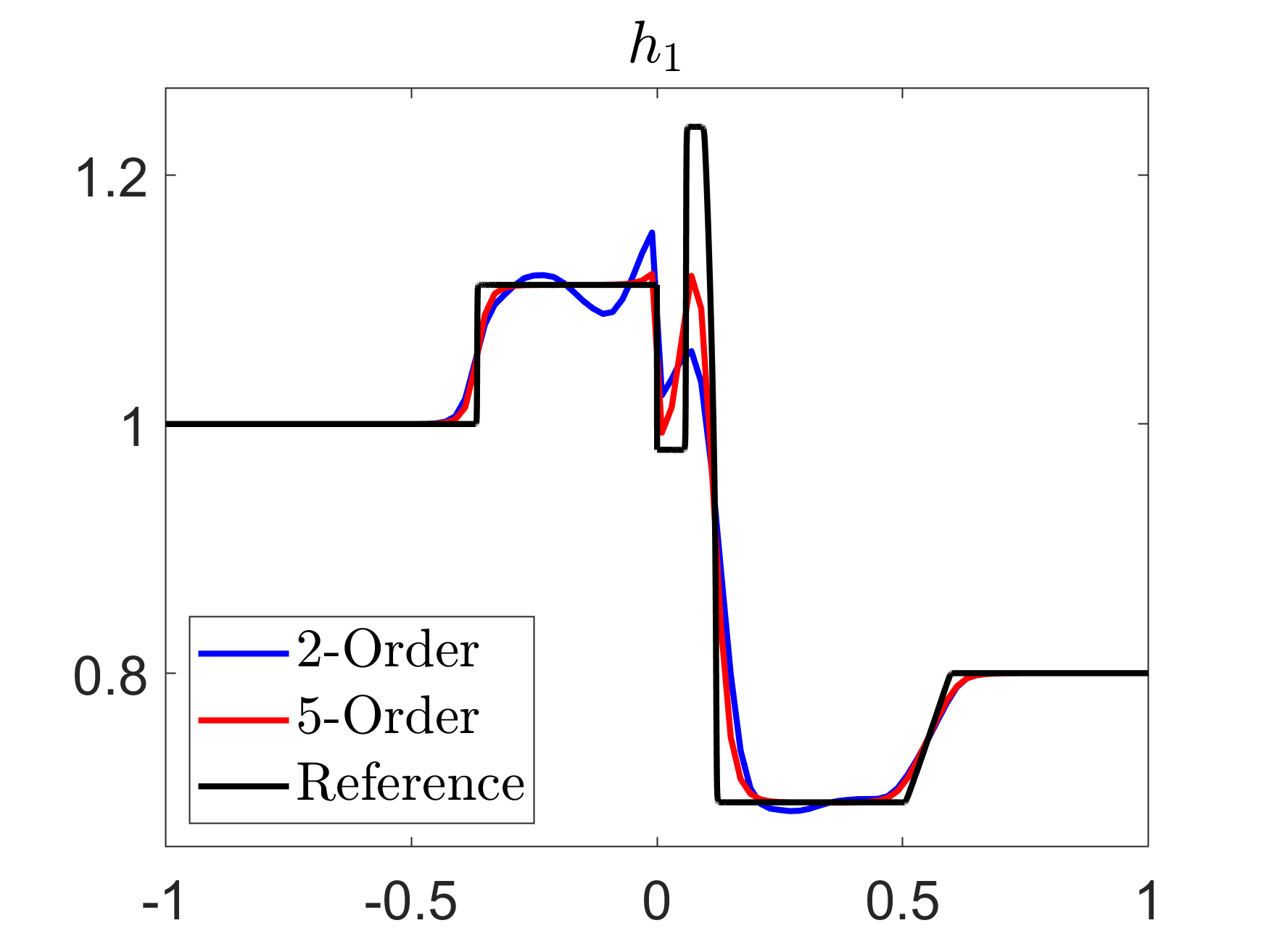}\hspace{0.5cm}
            \includegraphics[trim=0.8cm 0.3cm 1.cm 0.2cm, clip, width=5.cm]{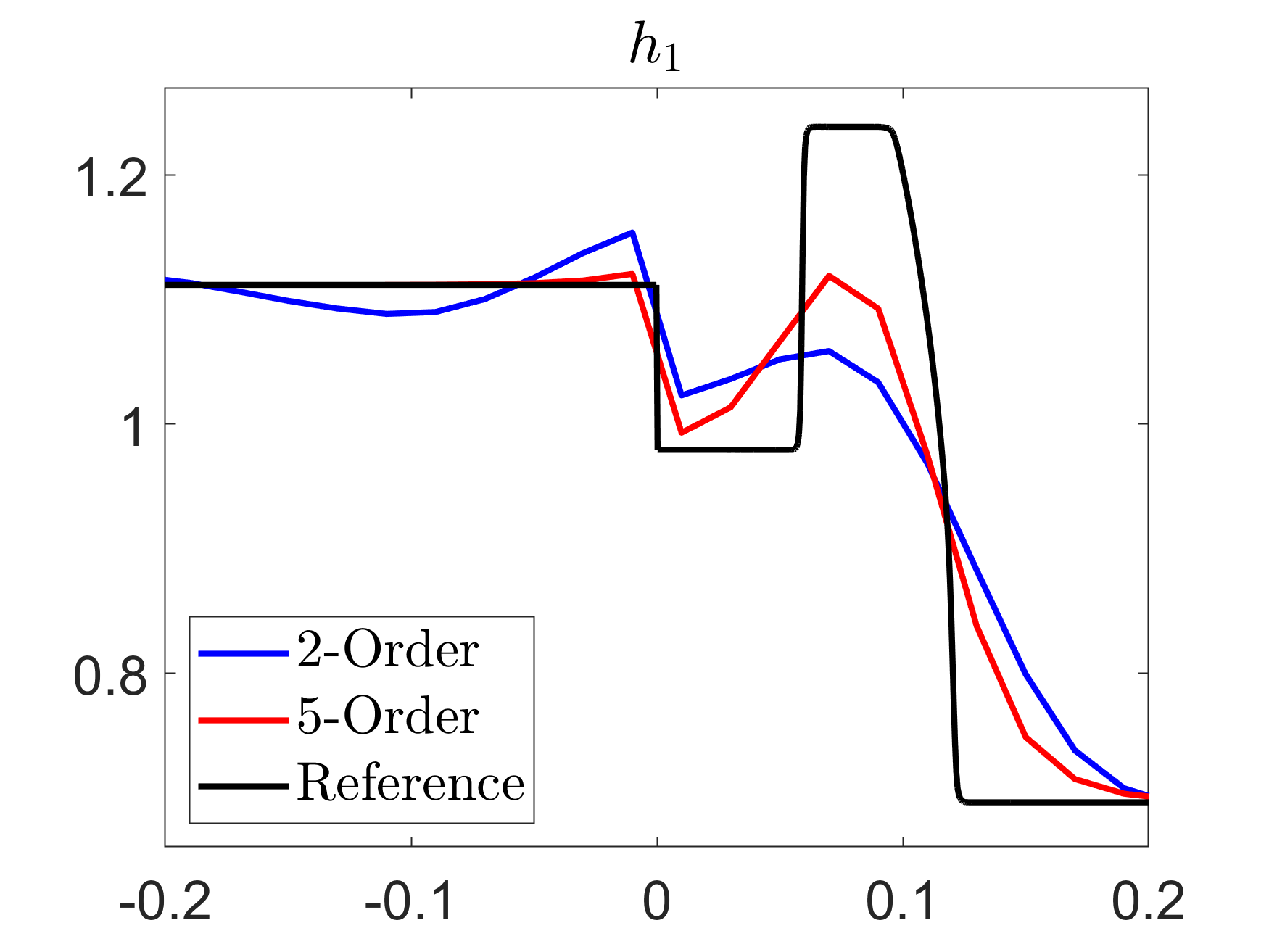}}
\vskip5pt
\centerline{\includegraphics[trim=0.8cm 0.3cm 1.cm 0.2cm, clip, width=5.cm]{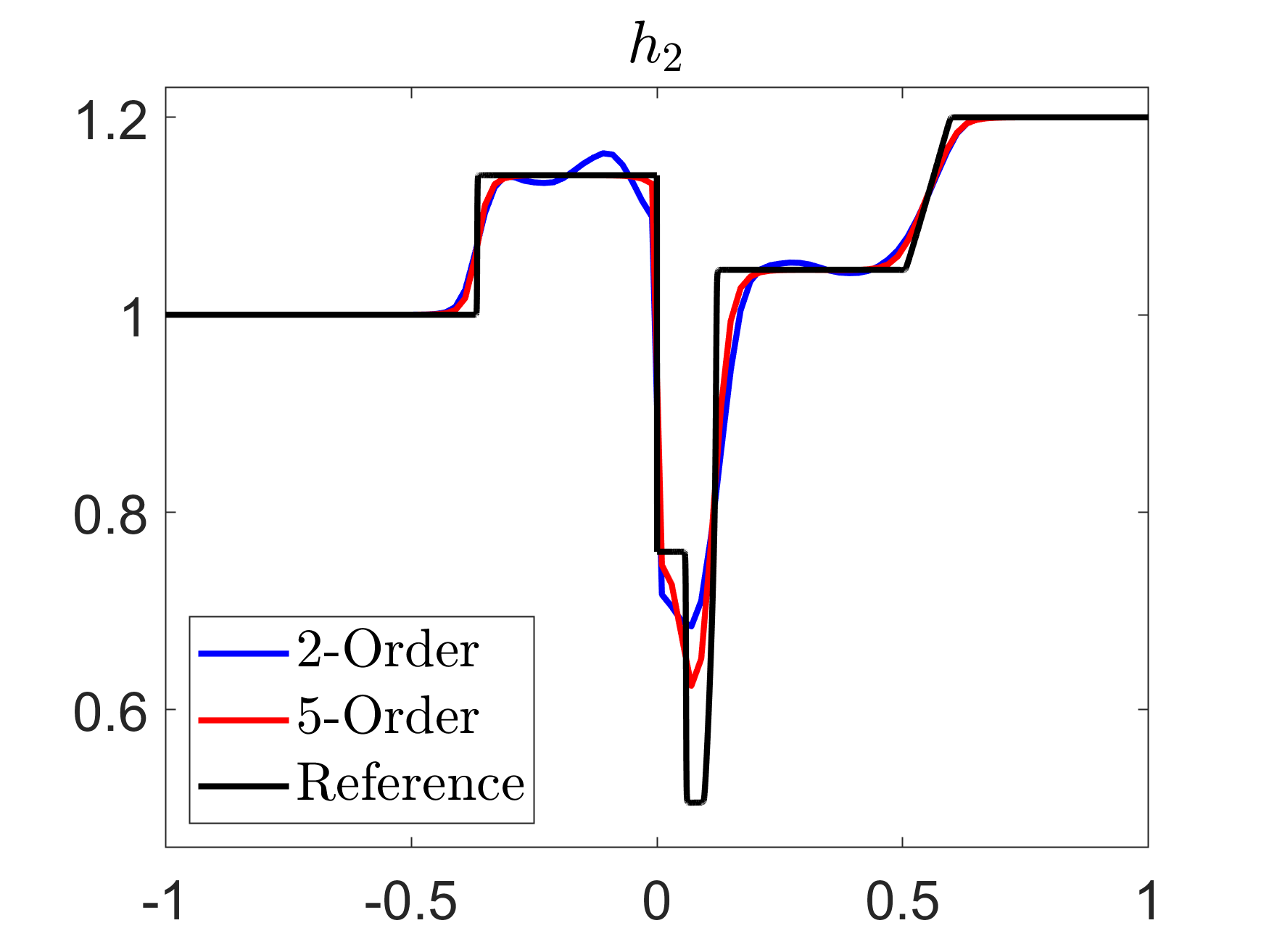}\hspace{0.5cm}
            \includegraphics[trim=0.8cm 0.3cm 1.cm 0.2cm, clip, width=5.cm]{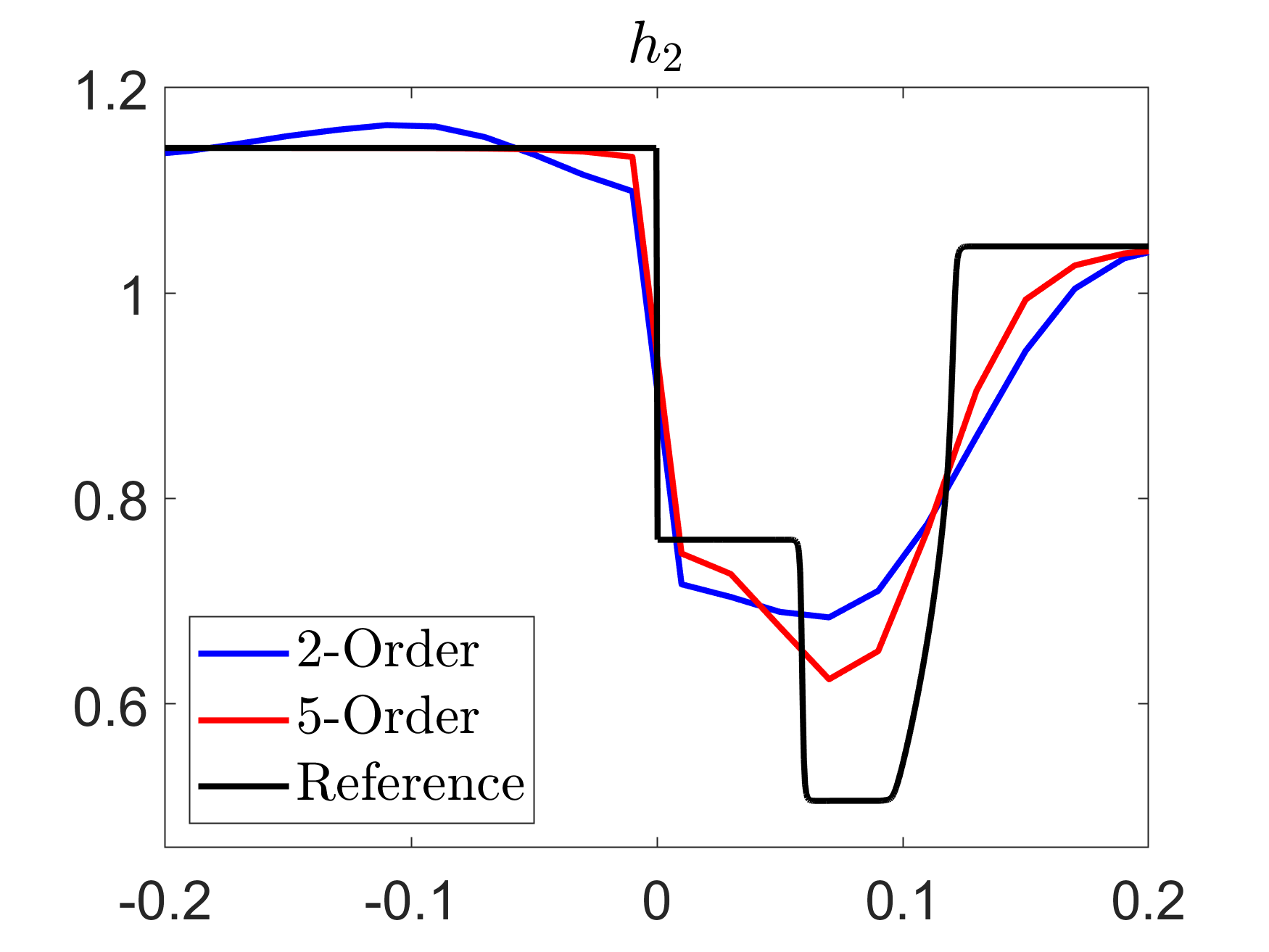}}
\caption{\sf Example 6 (Test 1): $h_1$ and $h_2$ (left column) and zoom at $x\in[-0.2,0.2]$ (right column).\label{fig12}}
\end{figure}

\begin{figure}[ht!]
\centerline{\includegraphics[trim=0.8cm 0.3cm 0.8cm 0.2cm, clip, width=5.cm]{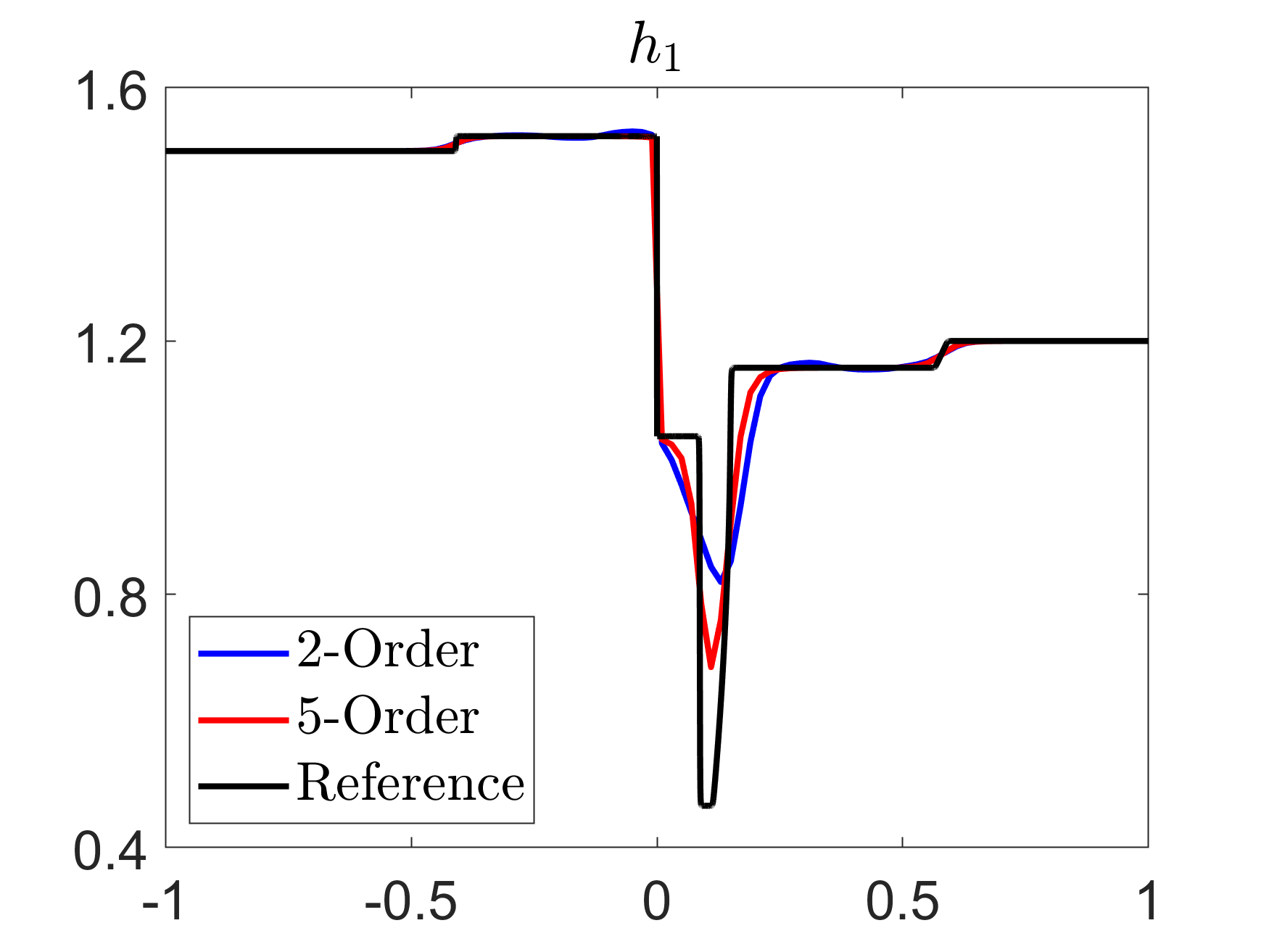}\hspace{1cm}
            \includegraphics[trim=0.8cm 0.3cm 0.8cm 0.2cm, clip, width=5.cm]{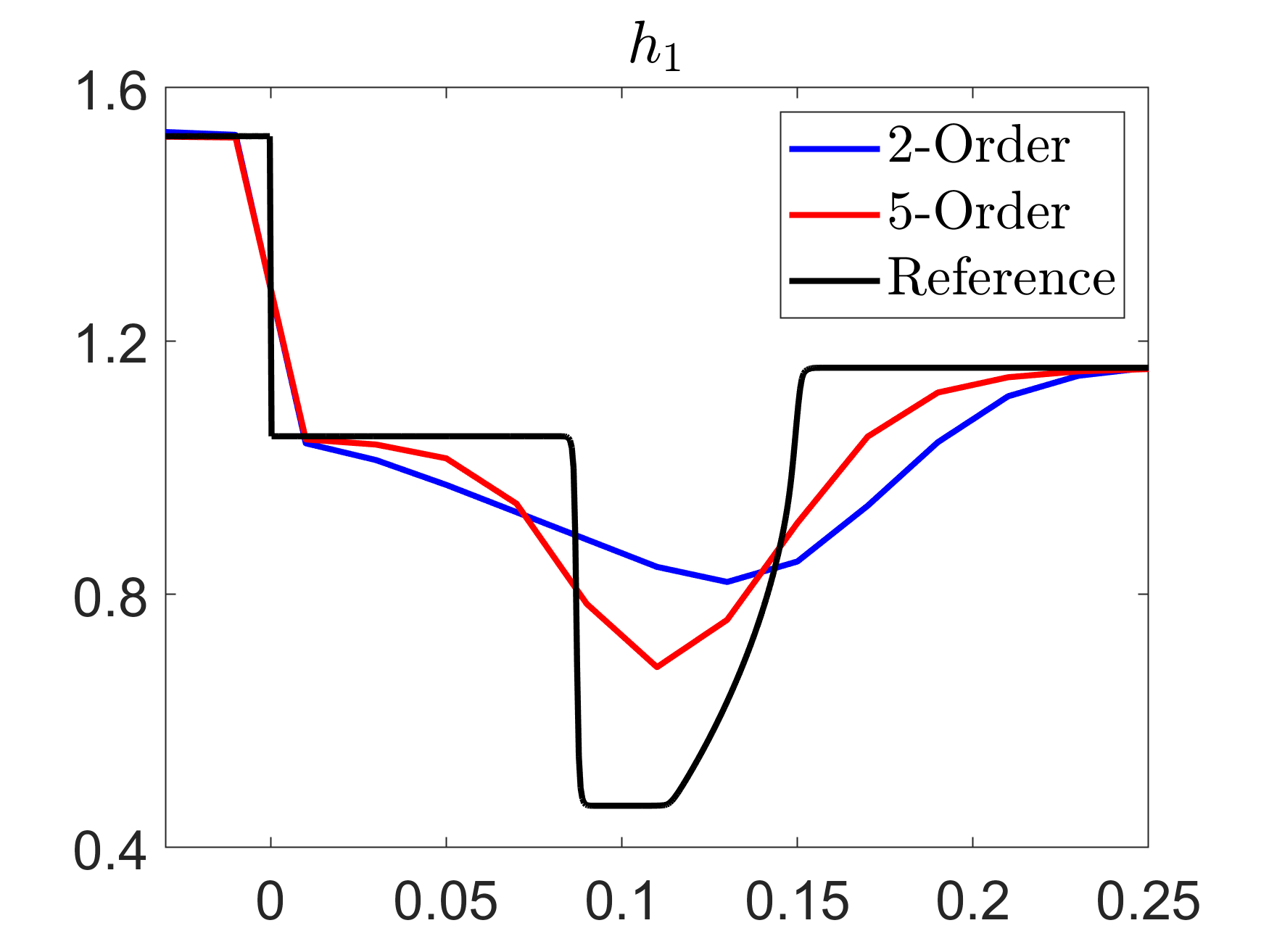}}
\vskip5pt
\centerline{\includegraphics[trim=0.8cm 0.3cm 0.8cm 0.2cm, clip, width=5.cm]{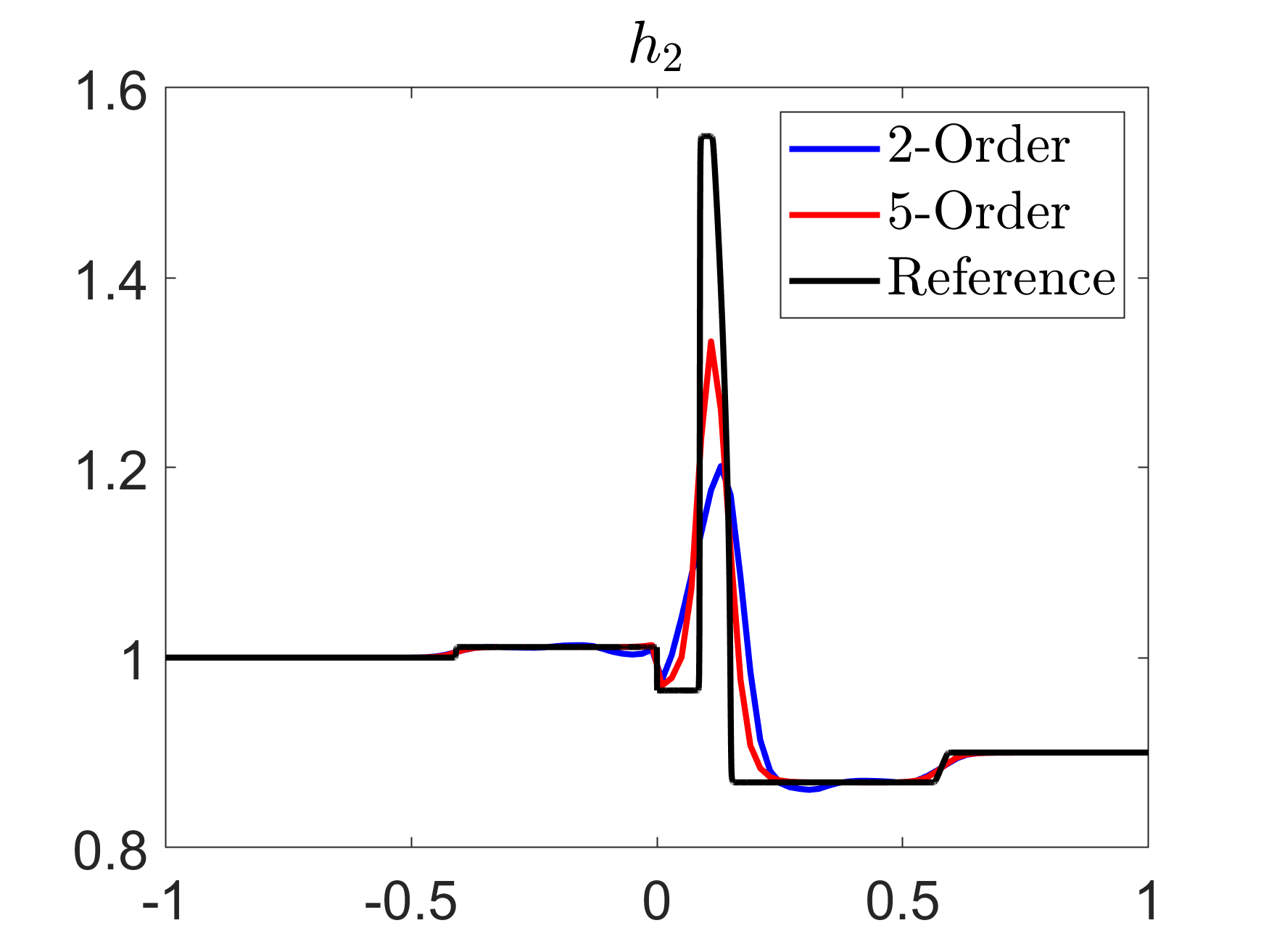}\hspace{1cm}
            \includegraphics[trim=0.8cm 0.3cm 0.8cm 0.2cm, clip, width=5.cm]{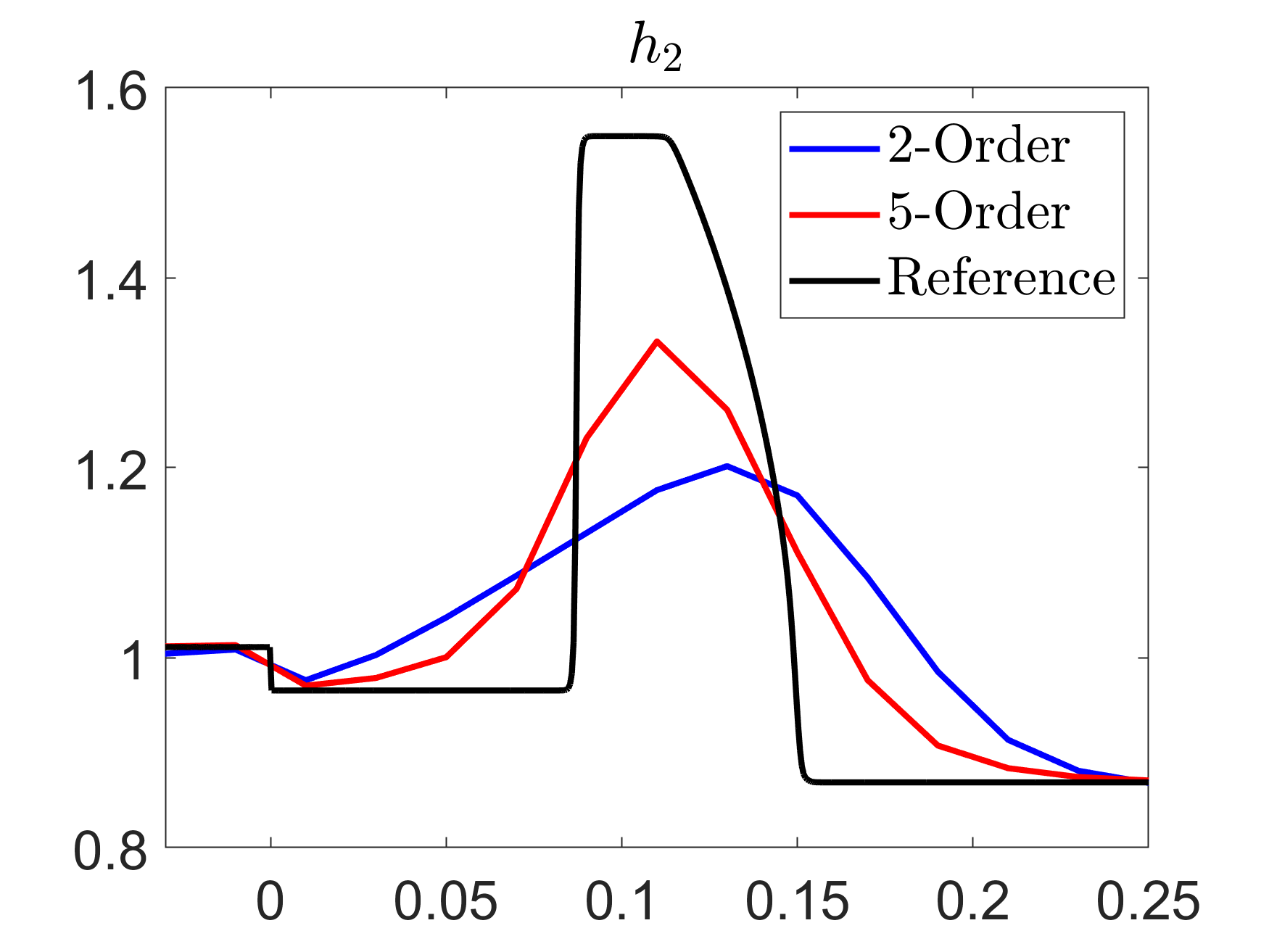}}
\caption{\sf Example 6 (Test 2): $h_1$ and $h_2$ (left column) and zoom at $x\in[-0.03,0.25]$ (right column).\label{fig13}}
\end{figure}

The numerical solutions computed by the 2-Order and 5-Order Schemes until time $t=0.1$ on a uniform mesh with $\dx=1/50$ are plotted in
Figures \ref{fig12} and \ref{fig13} together with the reference solution computed by the 5-Order Scheme on a much finer mesh with
$\dx=1/2000$. As one can see, the results obtained by the 5-Order Scheme are sharper and less oscillatory compared to those computed by the 
2-Order Scheme.

\section{Conclusions}\label{sec5}
In this paper, we have extended the second-order finite-volume flux globalization based well-balanced (WB) path-conservative central-upwind
(PCCU) schemes to fifth order of accuracy via the framework of the finite-difference alternative weighted essentially non-oscillatory
(A-WENO) schemes. The developed fifth-order flux globalization based WB A-WENO PCCU schemes have been applied to two nonconservative
systems---the nozzle flow system and two-layer shallow water equations. We have tested the fifth-order schemes on a number of numerical
examples and demonstrated that they clearly outperform their second-order counterparts.

\section*{Acknowledgements}
The work of S. Chu was supported in part by the DFG (German Research Foundation) through HE5386/19-3, 27-1. The work of A. Kurganov was supported in part by NSFC grant 12171226 and by the fund of the Guangdong Provincial Key Laboratory of Computational Science and Material Design (No. 2019B030301001).

\end{document}